\documentclass[preprint,12pt]{elsarticle}

\usepackage{amssymb,amsmath,amsthm,amsfonts}
\usepackage{breqn}
\usepackage{xcolor}
\newcommand{\Z}{\ensuremath{\mathbb{Z}}}

\newcommand{\C}{\ensuremath{\mathbb{C}}}

\newcommand{\R}{\ensuremath{\mathbb{R}}}
\newcommand{\Q}{\ensuremath{\mathbb{Q}}}
\usepackage{mathrsfs}
\usepackage{mathtools}
\usepackage{enumitem} 
\usepackage{hyperref}
\usepackage{graphicx}
\usepackage{subcaption}
\usepackage{etoolbox}
\apptocmd{\sloppy}{\hbadness 10000\relax}{}{}
\numberwithin{equation}{section}
\allowdisplaybreaks

\newtheorem{theorem}{Theorem}[section]
\newtheorem{lemma}{Lemma}[section]
\newtheorem{corollary}{Corollary}[section]

\journal{European Journal of Combinatorics}

\begin{document}

\begin{frontmatter}



\title{Full asymptotic expansion for orbit-summable quadrant walks and discrete polyharmonic functions}


\author[inst1,inst2]{Andreas Nessmann}

\affiliation[inst1]{organization={Institut für Diskrete Mathematik und Geometrie},%Department and Organization
            addressline={Technische Universität Wien},{Austria}}

\affiliation[inst2]{organization={Institut Denis Poisson},%Department and Organization
            addressline={Université de Tours et Université d'Orléans},{France} }
            \fntext[label1]{This project has received funding from the European Research Council (ERC) under the European Union's Horizon 2020 research and innovation programme under the Grant Agreement No. 759702.}

\begin{abstract}
Enumeration of walks with small steps in the quadrant has been a topic of great interest in combinatorics over the last few years. In this article, it is shown how to compute exact asymptotics of the number of such walks with fixed start- and endpoints for orbit-summable models with finite group, up to arbitrary precision. The resulting representation greatly resembles one conjectured for walks starting from the origin in~2020 by Chapon, Fusy and Raschel, differing only in terms appearing due to the periodicity of the model. We will see that the dependency on start- and endpoint is given by discrete polyharmonic functions, which are solutions of $\triangle^n v=0$ for a discretisation $\triangle$ of a Laplace-Beltrami operator. They can be decomposed into a sum of products of lower order polyharmonic functions of either the start- or the endpoint only, which leads to a partial extension of a theorem by Denisov and Wachtel.

\end{abstract}


\end{frontmatter}



The author would like to thank Kilian Raschel and Michael Drmota for a lot of valuable input and many helpful discussions.
\section{Introduction and Motivation}
Enumeration of lattice paths has by now become a standard problem in combinatorics. In particular the asymptotics of lattice path counting problem in the half- and quarter plane have been addressed in a variety of works. To do so, the by now standard approach as used in~\cite{BF02} is to consider the generating function 
\begin{equation}
    Q(x,y;t)=\sum_{n\in\mathbb{N}}\sum_{k,l\in \Z}x^ky^lt^n q(0,B;n),
\end{equation}
where $q(0,B;n)$ denotes the (possibly weighted) number of lattice paths starting at $(0,0)$ and ending at a point $B=(k,l)$. While for walks in the entire plane this contains negative powers of $x$ and $y$ and is thus not a generating function in these two variables, it is one in $t$. Utilizing a functional equation for $Q(x,y;t)$, it is then often possible to extract information about the asymptotic behaviour of $q(0,B;n)$, as was done for instance in~\cite{BF02}, where the authors gain various asymptotic expressions for the case of a directed model (i.e.~the first component of a step is always positive) closely tied to the zeros of the so-called kernel. In the undirected case, restricting ourselves to walks in the quarter plane, asymptotics of the coefficients of $Q(1,1;t)$ or $Q(0,0;t)$ (that is, of the walks starting at the origin and ending either anywhere or back at the origin) have been computed for some cases in~\cite{FR12}, using a complex boundary value problem, and for a family of models related to the Gouyou-Beauchamps model in~\cite{CMMR17} via Analytic Combinatorics in Several Variables (ACSV,~\cite{ACSV}). In~\cite{BCHKP}, the authors give hypergeometric expressions for the generating function in the $19$ cases where the group is finite (see~\cite{MBM}) and the model not algebraic, and show that it is transcendental. Additionally, they find explicit expressions for some specializations of the generating function for the $4$ algebraic models. To do so, they use orbit summation paired with creative telescoping. This goes in a similar direction as~\cite{BRS14}, and in general the question about the algebraic properties of the generating function $Q(x,y;t)$. This function is known to be D-finite in $x,y,t$ if and only if the group is finite, which was first conjectured in~\cite{MBM} and subsequently proven by various authors, for instance~\cite{BKvH10,MR09,KR12,Singer}.\\
There have also been extensions of the usual setting, where either steps are not required to be small or not to be homogeneous~\cite{FR17,KY15}, or where the setting is not in the quarter plane but some other domain, e.g. in $3$ or more dimensions~\cite{BBKM16,BHK21,MW19,F14,BW19}.\\
Instead of counting all walks of a certain length, however, one could also sort them by their endpoint. For some standard models, like for example the Simple Walk $\left(\mathcal{S}=\{\leftarrow,\uparrow,\rightarrow,\downarrow\}\right)$, or the Tandem Walk $\left(\mathcal{S}=\{\leftarrow,\uparrow,\searrow\}\right)$, explicit formulas for $q(0,B;n)$ are given in~\cite{MBM}, making use of coefficient extraction and orbit summation. One of the main results in the aforementioned article is that one can make use of orbit-summation in order to obtain an expression for the generating function $Q(x,y;t)$ in order to show it is D-finite, but it turns out that this approach can be utilized to compute asymptotics of $q(0,B;n)$ to arbitrary high order instead. For highly symmetric models, this is used in~\cite{MM16} in order to find diagonal representation of the generating functions, and then use ACSV to compute the asymptotics of the coefficients of $Q(1,1;t)$. In~\cite{MW19}, this is then refined to all orbit-summable models with finite group. Since the principal idea behind the saddle point method and ACSV is in fact rather similar (even though the latter is more general and uses much heavier machinery), it seems likely that their approach could be generalized to compute exact asymptotics of $q(k,l;n)$ as well. A general first order approximation with fixed start- and endpoint is given in~\cite{Denisov}, where using coupling with Brownian motion the authors show that, under some moment conditions, we have \begin{align}\label{eq:denisov}
    q(A,B;n)\sim \frac{V(A)\tilde{V}(B)}{n^c}+o\left(n^{-c}\right),
\end{align}
where $q(A,B;n)$ is the number of paths with starting- and endpoints $A=(u,v)$ and $B=(k,l)$ respectively, and $V(A), \tilde{V}(B)$ are discrete harmonic functions in $u,v$ and $k,l$ respectively. Note that these moment conditions are always satisfied for bounded step sets, which includes all cases considered in this article.\\
In~\cite{Poly}, the authors used methods from~\cite{BS97} to show that for a Brownian motion, the asymptotic expansion of the heat kernel with respect to time contains (continuous) polyharmonic functions. They noticed as well that a discrete analogue of this exists for the Simple Walk, and conjectured that this might also be the case for other models. Using orbit summation together with a saddle point method, we will show this to be true in a more general context in Thm.~\ref{thm:mainthm} and Thm.~\ref{thm:startingpoint}. In particular, for all orbit-summable models (see Section~\ref{sec:prelims}), we obtain for any given $m\in\mathbb{N}$ asymptotics of the form \begin{align}\label{eq:intro:asymptoticseries}
q(A,B;n)=\frac{\gamma^n}{n^c}\left[\sum_{p=1}^{m-1} \frac{v_p(A,B)\sum_{i=1}^r \alpha_i^{u-k}\beta_i^{v-l}\zeta_i^n}{n^p}+\mathcal{O}\left(\frac{1}{n^{m}}\right)\right],
\end{align}
where $c\in\mathbb{N}$, $\gamma\in\mathbb{R}$ and the $\alpha_i,\beta_i,\zeta_i$ are roots of unity. The coefficients $v_p(x,y)$ are so-called discrete polyharmonic functions (for a formal definition see Section~\ref{sec:prelims}) of order $p$ in $B$ and, in reversed direction, in $A$. If the model has no drift then they are polynomials, otherwise they contain additional exponential factors.  
Compared to~(\ref{eq:intro:asymptoticseries}), the expansion appearing in~\cite{Poly} is different in the lack of the $\alpha_i,\beta_i,\zeta_i$, which appear due to periodicity properties: at certain points we might have $q(A,B;n)=0$, and at those points some kind of cancellation needs to occur.\\
The structure of the polyharmonic functions $v_p(x,y)$ will be described more closely in Thm.~\ref{thm:startingpoint2}, where it will be shown that we can write \begin{align}\label{eq:intro_decomposition}
    v_p(A,B)=\sum_{i=1}^k h_{p,i}(A)g_{p,i}(B),
\end{align}
where the index $p$ is the same as in \eqref{eq:intro:asymptoticseries} and the $h_{i}(A)$ and $g_{i}(B)$ are (in the case of $y$ adjoint) polyharmonic of degree at most $p$ in $A$ and $B$ respectively (this will be detailed in Thm.~\ref{thm:startingpoint2}), and the number of summands $k$ depends on $p$. For $p=1$ this reduces precisely to~(\ref{eq:denisov}), so in a certain sense one could view Thm.~\ref{thm:startingpoint2} as an extension (albeit under much stronger conditions) of~(\cite[Thm.~6]{Denisov}).\\
While the results of this article are stated for orbit-summable small step models in two dimensions, they do generalize to higher dimensions and -- with a slight condition -- to orbit-summable large step models as well, see the remarks after Thm.~\ref{thm:mainthm} as well as \ref{app:largesteps} and \ref{app:3dim}.\\
This article is structured as follows: \begin{enumerate}
\item In Section~\ref{sec:prelims}, a short overview over some definitions and the tools utilized will be given.
\item In Section~\ref{sec:quadrantwalks_origin}, it will be shown that models with finite group allowing for orbit summation in a manner similar as in~\cite{MBM} satisfy an asymptotic relation of the form (\ref{eq:intro:asymptoticseries}), and the example of the Gouyou-Beauchamps model will be worked out in detail. In particular it will be explicitly shown how to compute the constants $\gamma$, $c$, $\alpha_i,\beta_i,\zeta_i$ and the functions $v_{p}(B)$ in (\ref{eq:intro:asymptoticseries}).
\item In Section~\ref{sec:quadrantwalks_general}, we will consider the same problem where instead of $(0,0)$, we start our paths at an arbitrary point $(u,v)$. It turns out, maybe not surprisingly, that due to the symmetry of the problem the resulting solution is similar to that obtained in Section~\ref{sec:quadrantwalks_origin}, leading to Thm.~\ref{thm:startingpoint}. Continuing from there, we can then find a decomposition of the $v_p$ as in (\ref{eq:intro_decomposition}) in Thm.~\ref{thm:startingpoint2}. The resulting decomposition of the asymptotic terms for the Simple Walk is given in Section~\ref{sec:startingpoints:SW}, using an explicit basis of polyharmonic functions constructed in~\cite{Nes22}.
\item In~\ref{app:largesteps} and~\ref{app:3dim}, examples of the method applied to a model with large steps and a three-dimensional model are given. Additional examples of the decomposition of the polyharmonic coefficients as in Section~\ref{sec:startingpoints:SW} for the Gouyou-Beauchamps model and the Tandem Walk as well as the first three terms of the asymptotics for all the $19$ unweighted, orbit-summable models are given in the preprint of this article~\cite[App.~C,~D]{Nes23}.

\end{enumerate}
\section{Preliminaries}\label{sec:prelims}
\subsection{Walks in the Quarter Plane}
In order to count lattice paths in the quarter plane, first of all we need a set $\mathcal{S}\subseteq\mathbb{Z}^2$ of permissible steps, together with a family of weights $\left(\omega_s\right)_{s\in\mathcal{S}}$. A lattice path of length $n$ is then a sequence of points $(x_0,\dots, x_n)\subseteq \mathcal{Q}^n$, with $\mathcal{Q}:=\mathbb{Z}_{\geq 0}\times\mathbb{Z}_{\geq 0}$ being the quarter plane, such that $x_k-x_{k-1}\in\mathcal{S}$ for all $k\in\{1,\dots,n\}$. We will count such a path by its weight, that is, the product $\prod_{k=1}^n\omega_{x_k-x_{k-1}}$.\\
We will in the following assume that: \begin{enumerate}
\item our step set consists of small steps only, i.e.~$\mathcal{S}\subseteq\{-1,0,1\}^2\setminus\{(0,0)\}$,
\item our step set is non-degenerate, i.e.~there is no (possibly rotated) half-plane containing all allowed steps.
\end{enumerate}
In order to keep the notation short, we will in the following sometimes denote by $\mathcal{S}$ not only the set of allowed steps, but also the associated weights $\left(\omega_s\right)_{s\in\mathcal{S}}$.\newline 
Let $q((u,v),(k,l);n)$ be the (weighted) number of paths of length $n$ from $(u,v)$ to the endpoint $(k,l)$. It can then be shown (see e.g.~\cite[Lemma~4]{MBM}) that the generating function \begin{align}
Q_{u,v}(x,y;t)=\sum_{n\geq 0}t^n\sum_{k,l\in\mathbb{N}}x^ky^l q((u,v),(k,l);n)\end{align}
satisfies the functional equation \begin{small}\begin{multline}\label{eq:prelims:feq}
K(x,y;t)Q_{u,v}(x,y;t)\\=x^{u+1}y^{v+1}+\underbrace{K(x,0;t)Q_{u,v}(x,0;t)}_{=:A(x;t)}+\underbrace{K(0,y;t)Q_{u,v}(0,y;t)}_{=:B(y;t)}-\underbrace{K(0,0;t)Q_{u,v}(0,0;t)}_{=:C(t)},
\end{multline}
\end{small}
where \begin{align}
K(x,y)=xy\left[1-tS(x,y)\right],\end{align}
where $S(x,y)$ is the step-counting Laurent polynomial \begin{align}
S(x,y)=\sum_{(i,j)\in\mathcal{S}}\omega_{(i,j)}x^iy^j.
\end{align}
For ease of notation, if $(u,v)=(0,0)$, then we will just write $q(k,l;n)$ and $Q(x,y;t)$.\\
Given that we consider paths in the quarter plane only, we see that the series $Q_{u,v}(x,y;t)$ are indeed power series in $x,y$ and $t$; in particular they are power series in $t$ with polynomial coefficients in $x,y$.\\
Given some power series $A(x,y)=\sum_{i,j\in\mathbb{N}}a_{i,j}x^iy^j$, we denote by $[x^n]$ the linear operator extracting the $n$-th coefficient, i.e.~$[x^n]A(x,y)=\sum_{j\in\mathbb{N}}a_{n,j}y^j$, and $[x^ny^m]A(x,y)=a_{n,m}$. In the same manner, we define for Laurent series the operator $[x^>]:=\sum_{n>0}^\infty x^n[x^n]$ extracting all positive powers with their coefficients, and the operator $[x^\geq]$ extracting all nonnegative ones.\newline
As we allow small steps only and our step set is non-degenerate, we can write the kernel $K(x,y)$ as \begin{align}\label{eq:kernelform1}
    K(x,y)&=a(x)y^2+b(x)y+c(x)\\\label{eq:kernelform2}
    &= \tilde{a}(y)x^2+\tilde{b}(y)x+\tilde{c}(y),
\end{align}
with $a(x),b(x),c(x),\tilde{a}(y),\tilde{b}(y),\tilde{c}(y)$ all being non-zero. Consider now the two birational transformations \begin{align}
\label{eq:group1}
\Phi: (x,y)\mapsto &\left(x,y^{-1}\frac{c(x)}{a(x)}\right),\\
\label{eq:group2}
\Psi: (x,y)\mapsto &\left(x^{-1}\frac{\tilde{c}(y)}{\tilde{a}(y)},y\right).
\end{align}
These two transformations, which clearly depend on our step set $\mathcal{S}$, generate the so-called \textbf{group of a model}, denoted by $\mathcal{G}$. This group can be either finite or infinite. If $\mathcal{G}$ is finite, then as both $\Phi,\Psi$ are involutions, any element $g\in\mathcal{G}$ can be written as either $\left(\Phi\circ\Psi\right)^k$, or as $\Psi\circ\left(\Phi\circ\Psi\right)^k$ for some $k$. We define $\operatorname{sgn}(g)=1$ in the first, and $\operatorname{sgn}(g)=-1$ in the second case. The study of this group has been central to many results about lattice walks in the quarter plane, see e.g.~\cite{MBM,Singer,SingerGenus0,FR12}. The main reason for this is that if we let
\begin{align}
k(x,y):=\frac{K(x,y)}{xy}=1-tS(x,y),
\end{align}
then it follows that $k(x,y)$ is invariant under $\mathcal{G}$. This observation is the starting point for orbit summation methods as is done in~\cite{MBM}: we rewrite (\ref{eq:prelims:feq}) (leaving out the dependency on $t$ for easier readability) as \begin{align}
    xyk(x,y)Q(x,y)=xy+A(x)+B(y)-C.
\end{align}
Assuming now that the group is finite, and picking any $g\in\mathcal{G}$, we hence obtain \begin{align}\label{eq:prelims:osum1}
    g(xy)k(x,y)Q(g(x),g(y))=g(xy)+A(g(x))+B(g(y))-C.
\end{align}
Multiplying (\ref{eq:prelims:osum1}) with $\operatorname{sgn}(g)$ and taking the sum over all elements $g\in\mathcal{G}$, all terms $A(g(x))$ and $B(g(y))$ cancel (note that both $\Psi,\Phi$ change only one variable each, while $\operatorname{sgn}$ switches sign), we obtain (see~\cite[Prop.~5]{MBM}) \begin{align}
\sum_{g\in\mathcal{G}}\operatorname{sgn}(g) g(xy)Q(g(x),g(y))=\frac{1}{k(x,y)}\sum_{g\in\mathcal{G}}\operatorname{sgn}(g) g(xy).
\end{align}
It often turns out that often none of the terms of the form $Q(g(x),g(y))$ will contribute positive powers of both $x$ and $y$ except if $g=\operatorname{id}$. In this case, we say that the model is \textbf{orbit-summable}, and we can compute $Q(x,y;t)$ by coefficient extraction via \begin{align}\label{eq:prelims:osummable}
xyQ(x,y;t)=\left[x^>\right]\left[y^>\right]\frac{\sum_{g\in\mathcal{G}}\operatorname{sgn}(g)g(xy)}{k(x,y;t)}.
\end{align}
By~\cite[Lemma~2]{MBM}, we know that if a model $\mathcal{S}$ (consisting of both steps and weights $\omega_s$) has a finite group, then the model with reversed steps $\tilde{S}:=-\mathcal{S}$ has a finite group as well. Furthermore, by essentially the same argument, one can see that orbit-summability is retained when reversing the steps as well:
\begin{lemma}\label{lemma:orbitsummability_reversed}
    Given an orbit-summable model $\mathcal{S}$, the reversed model $\tilde{\mathcal{S}}$ is also orbit-summable.
\end{lemma}

\begin{proof}
    We note as in~\cite[Lemma~2]{MBM} that if $\Phi,\Psi$ are the generators of the group of $\mathcal{S}$, then the generators of the group of $\tilde{\mathcal{S}}$ are given by $\iota \circ\Psi\circ\iota,\iota\circ\Phi\circ\iota$, where $\iota(x,y)=\left(x^{-1},y^{-1}\right)$. Given any element $g$ in the group $\mathcal{G}$ of $\mathcal{S}$, denote by $\tilde{g}$ the corresponding element of the group $\tilde{\mathcal{G}}$ of $\tilde{\mathcal{S}}$. As $\iota$ is an involution, we see that $\tilde{g}=\iota\circ g \circ \iota$ for any such $g$, and in particular that $\tilde{g}(xy)=\iota\circ g\circ \iota (xy)$. Noticing that $g(xy)\in\mathbb{C}(x,y)$ we can deduce that $\tilde{g}(xy)=g(xy)$. From there, we obtain an equivalent equation to (\ref{eq:prelims:osummable}) by the very same combinatorial arguments.
\end{proof}
In the following, it will be convenient to consider models with zero drift, that is, where \begin{align*}\sum_{(i,j)\in\mathcal{S}}i\omega_{i,j}=\sum_{(i,j)\in\mathcal{S}}j\omega_{i,j}=0.\end{align*} If this is not the case, then we can utilize the \textbf{Cramer-transformation}: we multiply each weight $\omega_{i,j}$ by a factor of $\alpha^i\beta^j$, where we choose $\alpha,\beta$ such that the drift will be $0$. The existence of such an $\alpha,\beta$ is ensured for non-singular models, see e.g.~\cite[1.5]{Denisov}. The reason why this substitution is very convenient combinatorially is fairly simple; given the number $q((i,j),(k,l);n)$ of paths from $(i,j)$ to $(k,l)$ with $n$ steps weighted by old weights $\omega_{i,j}$, then for the equivalent $\hat{q}$ using the new weights we have\begin{align}
    \hat{q}((i,j),(k,l);n)=\alpha^{k-i}\beta^{j-l}q((i,j),(k,l);n).
\end{align}
Additionally, it turns out that the group of the model is, in a certain sense, invariant under this transformation, and in particular orbit-summability is preserved.
\begin{lemma}\label{lemma:cramer}
    Let $\mathcal{S}$ be a model and $\hat{\mathcal{S}}$ be a Cramer-transform of $\mathcal{S}$, with weights $\alpha,\beta$. Let $\mathcal{G},\hat{\mathcal{G}}$ be the respective groups with generators $\Phi,\Psi$ and $\hat{\Phi},\hat{\Psi}$. Lastly, let $\iota$ be the mapping $(x,y)\mapsto (\alpha x,\beta y)$. Then we have:\begin{enumerate}
        \item $\hat{\Phi}(x,y)=\iota^{-1}\circ \Phi \circ \iota$, $\hat{\Psi}(x,y)=\iota^{-1}\circ\Psi\circ\iota$; in particular $\mathcal{G}$ and $\hat{\mathcal{G}}$ are isomorphic, 
        \item $\hat{\mathcal{S}}$ is orbit-summable if and only if $\mathcal{S}$ is orbit-summable.
    \end{enumerate}
\end{lemma}
\begin{proof} We show the first part for $\Phi$ only, the statement for $\Psi$ will then follow by symmetry. We can define $\hat{c}(x)$ and $\hat{a}(x)$ as in (\ref{eq:kernelform1}). One finds that, by definition of $\hat{\mathcal{S}}$, we will have $\hat{a}(x)=\beta a(\alpha x)$ and $\hat{c}(x)=\beta^{-1}c(\alpha x)$. From this it follows that \begin{align}
    \iota^{-1}\circ\Phi\circ\iota (x,y)=\iota^{-1}\circ\Phi(\alpha x,\beta y)=\iota^{-1}\left[\alpha x,(\beta y)^{-1}\frac{c(\alpha x)}{a(\alpha x)}\right]\\
    =\iota^{-1}\left(\alpha x, \beta y^{-1}\frac{\hat{c}(x)}{\hat{a}(x)}\right)=\left( x, y^{-1}\frac{\hat{c}(x)}{\hat{a}(x)}\right)=\hat{\Phi}(x,y).
\end{align}
The isomorphism of $\mathcal{G}$ and $\hat{\mathcal{G}}$ follows immediately and is given by $g\mapsto \iota^{-1} \circ g \circ \iota \in\hat{\mathcal{G}}$. By this isomorphism, we also see immediately that positive and negative powers of $x,y$ are preserved, which is in turn all that matters for orbit-summability, hence we are done.
\end{proof}
Given a model with zero drift, one can consider the correlation coefficient $\theta$ (see~\cite{FR12,Conformal}) given by \begin{align}\label{eq:deftheta}
    \theta=\arccos\left(-\frac{\sum_{(i,j)\in\mathcal{S}} ij\omega_{i,j}}{\sqrt{\sum_{(i,j)\in\mathcal{S}}i^2\omega_{i,j}\cdot\sum_{(i,j)\in\mathcal{S}}j^2\omega_{i,j}}}\right).
\end{align}
It turns out that this correlation coefficient is closely linked to the properties of the group (its restriction to the surface $\mathcal{C}:=\{(x,y): K(x,y)=0\}$ is finite if and only if $\pi/\theta\in\Q$~\cite{FR11}), and also has a geometric interpretation (see~\cite{Conformal}). Additionally, in~\cite{Kilian_JEMS},~\cite{Hung}, it is shown that $\pi/\theta$ is tied to the growth of the positive harmonic function. Given a model with non-zero drift, we can first apply a Cramer-transform to get rid of the drift, and then define $\theta$ as above. 
\subsection{Discrete Polyharmonic Functions}
Given a step set $\mathcal{S}$ and a discrete function $f$ defined on the quarter plane $\mathcal{Q}$, we can define the Markov operator \begin{align}
    Pf(x):=\sum_{s\in\mathcal{S}}\omega_s f(x-s).
\end{align}
If we take the discrete random walk $(X_n)$ with the transition probabilities given by the reverse of $\mathcal{S}$, and the induced Markov chain $M_n:=f(X_n)$, then we can interpret the operator $P$ as the expectation $\mathbb{E}\left[M_{n+1}\mid M_n\right]$. One can then proceed to look at the expected change during a time step, weighted by a parameter $t$, which is given by \begin{align}
    \triangle f(x):=\left(P-t\operatorname{id}\right)f(x).
\end{align}
The operator $\triangle$ can be viewed as the discrete equivalent of a Laplace-Beltrami operator. We call a function \textbf{(discrete) $t$--harmonic} if \begin{enumerate}
    \item $\triangle f(x)=0$ for all $x\in\mathcal{Q}$,
    \item $f(x)=0$ for all $x\in\mathcal{Q}^c$, 
\end{enumerate}
where $\mathcal{Q}^c$ is the complement of $\mathcal{Q}$. Similarly, we call a function \textbf{(discrete) $t$--polyharmonic of degree $p$} if \begin{enumerate}
    \item \label{def:PHF_cond1}$\triangle^p f(x)=0$ for all $x\in\mathcal{Q}$,   
    \item \label{def:PHF_cond2} $f(x)=0$ for all $x\in\mathcal{Q}^c$. 
\end{enumerate}
If $t=1$, then we simply speak of harmonic and polyharmonic functions respectively. Note that strictly speaking, the polyharmonicity is given by the first conditions, while the second ones are Dirichlet boundary conditions. Due to the underlying combinatorics, however, one immediately sees that we are only interested in polyharmonic functions satisfying the Dirichlet problem as above. In particular ($1$-)harmonic and biharmonic functions have a variety of applications in the theory of stochastic processes and physics, see e.g.~\cite{Elasticity,HFApplications,CG93}. The occurrence of $t$-polyharmonic functions in the asymptotics of path counting problems, however, was noted only fairly recently in~\cite{Poly}. In~\cite{Nes22}, it is shown that in the zero drift case, the space of $1$-polyharmonic functions of order $n$ is isomorphic to $\mathbb{C}[[z]]^n$, and a basis consisting of polyharmonic functions with rational generating function is given.\newline
Given a model $\mathcal{S}$, we can also define a discrete Laplacian $\tilde{\triangle}$ for the model with directions reversed $\tilde{\mathcal{S}}$. Since $\tilde{\triangle}$ is the adjoint operator to $\triangle$ on the space $\mathcal{L}^2\left(\Z^2\right)$, we will call it the adjoint Laplacian. In Section~\ref{sec:startingpoints}, we will encounter functions of the form $f(x,y)$, for $x,y\in\Z^2$, which are polyharmonic in $x$ and adjoint polyharmonic in $y$. We will call such a function $f(x,y)$ \textbf{multivariate polyharmonic} of order $p$ if \begin{align*}
    \triangle^k\left(\tilde{\triangle}^{p-k}f\right)=0
\end{align*}
for all $0\leq k\leq p$. Note that one can verify by an easy computation using linearity that the ordering of the Laplacians does not matter as they commute, i.e.~we have \begin{align}
    \triangle\left(\tilde{\triangle}f(x,y)\right)=\tilde{\triangle}\Bigl(\triangle f(x,y)\Bigr).
\end{align} \newline
Lastly, discrete polyharmonic functions behave well with respect to Cramer transformations: via a short computation, one can show
\begin{lemma}\label{:cramer_PHF}
Let $\hat{\mathcal{S}}$ be a Cramer-transform of $\mathcal{S}$, with $\hat{\omega}_{i,j}=\alpha^i\beta^j\omega_{i,j}$. Let $\hat{\triangle}, \triangle$ be the associated Laplacians. We then have \begin{align}
    \hat{\triangle}\left[\alpha^{-k}\beta^{-l} f(k,l)\right]=\alpha^{-k}\beta^{-l}\triangle\left[f(k,l)\right].
\end{align}
\end{lemma}
This directly implies that we have a bijection between the polyharmonic functions w.r.t.~$\triangle$ and those w.r.t.~$\hat{\triangle}$, given by adding a factor of $\alpha^{-k}\beta^{-l}$.

\subsection{Saddle Points}
In the following we will use the saddle point method as described for instance in~\cite[VIII]{Flajolet}. In particular, we will be interested in saddle points of $S(x,y)$, which is a Laurent polynomial with only positive coefficients. We call a point $(x_0,y_0)$ a \textbf{dominant saddle point} if $(x_0,y_0)$ is a minimizer of $S(x,y)$ in $\mathbb{R}^+\times\mathbb{R}^+$. Note that this minimization property directly implies the usual defining property of a saddle point, namely that the orthogonal directional derivatives vanish. Additionally, note that $S(x_0,y_0)>0$, due to the positivity of coefficients.
\begin{lemma}
    For any non-degenerate model, a dominant saddle point exists.
\end{lemma}
\begin{proof}
    Since the model is non-degenerate, we know that $S(x,y)$ is coercive on $\mathbb{R}^+\times\mathbb{R}^+$ (it goes to infinity wherever we approach the boundary; see also~\cite[1.5]{Denisov}). Therefore, $S(x,y)$ must attain its minimum at some point $(x_0,y_0)$, which is then by definition a dominant saddle point.
\end{proof}

By a short computation, one can check that \begin{lemma}\label{lemma:saddlepoint_driftzero}
The dominant saddle point $s_0=(x_0,y_0)$ is equal to $(1,1)$ if and only if the drift of the model is $0$.
\end{lemma}
When using a Cauchy type integral to compute asymptotics later on, we will integrate over the torus $\{(x,y): |x|=x_0,|y|=y_0\}$. We will call a point $(x_i,y_i)$ on this torus a \textbf{saddle point associated with $(x_0,y_0)$}, if $|S(x_i,y_i)|=S(x_0,y_0)$.\newline
Note that by the positivity of the coefficients of $S(x,y)$, we know that $S(x,y)$ takes the maximum (absolute) value on the aforementioned torus at $(x_0,y_0)$. Indeed, for any other point $(x_i,y_i)$ we can have $S(x_i,y_i)=S(x_0,y_0)$ only if there is some $\zeta_i\in\mathbb{C}, |\zeta_i|=1$ such that each monomial term of $S(x_i,y_i)$ differs from the corresponding term of $S(x_0,y_0)$ by this factor $\zeta_i$.\newline
Choosing $\alpha_i,\beta_i$ such that $(x_i,y_i)=(\alpha_i x_0,\beta_i y_0)$, we must therefore have $\alpha_i^k\beta_i^l=\zeta_i$ for all $(k,l)$ such that $\omega_{k,l}\neq 0$. From here, it is not difficult to see that $\alpha_i,\beta_i$ and $\zeta_i$ must be roots of unity. It is also clear that there can only be finitely many such $\zeta_i$, in a one-to-one correspondence with finitely many pairs $(x_i,y_i), 0\leq i\leq l$ (the maximum $l$ appearing for the $19$ unweighted orbit-summable models is $3$, see~\cite[App.~D]{Nes23}). \newline
One can check directly that the $(x_i,y_i)$ are indeed saddle points of $S(x,y)$, moreover, as we will see in Lemma~\ref{lemma:periodicity}, the local behaviour of $S(x_0,y_0)$ and $S(x_i,y_i)$ is the same up to the factor $\zeta_i$. By the same reasoning as in~\cite[VIII]{Flajolet}, it turns out that when computing the asymptotics of the coefficients via the Cauchy formula, the main contributions to the contour integral come from the points $(x_i,y_i)$, as the modulus of $S(x,y)$ will be smaller elsewhere, leading to exponentially smaller terms.\newline
Note that while this article considers mainly the $2$-dimensional case, all the definitions above extend to more dimensions in a natural manner; only the structure of the group becomes more complicated as we will have more than two transformations, see e.g.~\cite{BBKM16}. 

\section{Quadrant walks starting at the origin}
\label{sec:quadrantwalks_origin}
\subsection{Full asymptotic expansion}\label{sec:maintheorem}
The goal of this section is to compute the asymptotics of orbit-summable lattice walks from the origin to an arbitrary but fixed point in the quarter plane, and in particular to show the following:
\begin{theorem}\label{thm:mainthm}
Let $\mathcal{S}$ be a step set satisfying the general assumptions stated in Section~\ref{sec:prelims} and such that $\mathcal{S}$ is orbit-summable, i.e.~\begin{align*}
    xyQ(x,y;t)=\left[x^>\right]\left[y^>\right]\frac{N(x,y)}{k(x,y;t)},
\end{align*}
where \begin{align*}
    N(x,y):=\sum_{g\in\mathcal{G}}\operatorname{sgn}(g)g(xy).
\end{align*}
Suppose that $s_0=(x_0,y_0)$ is a dominant saddle point, with associated other saddle points $s_i$, $1\leq i\leq r$ (meaning that we consider $r+1$ saddle points in total). Furthermore, let $\alpha_i,\beta_i,\zeta_i$ be the roots of unity as constructed in Section~\ref{sec:prelims}. Then, there is a constant $c\in\mathbb{N}$, a constant $\gamma>0$, and $\gamma$-polyharmonic functions $v_p$ of degree $p$ such that for any $m\in\mathbb{N}$ we have \begin{align}\label{eq:mainthm:asymptotics}
q(0,(k,l);n)=\frac{\gamma^n}{n^c}\left[\sum_{p=1}^{m-1}\frac{v_p(k,l)\sum_{i=1}^r \alpha_i^k\beta_i^l\zeta_i^n}{n^p}+\mathcal{O}\left(\frac{1}{n^m}\right)\right].
\end{align}
The polyharmonic functions $v_p(k,l)$ are polynomials precisely if the drift of the model is zero, else they contain an additional factor of $x_0^{-k}y_0^{-l}$, with $(x_0,y_0)$ the dominant saddle point. They can be computed explicitly via a Cauchy-type integral.\\
Lastly, the constant $c$ can be expressed using the correlation coefficient $\theta$ defined in (\ref{eq:deftheta}) via $c=\pi/\theta$, and for the exponential growth $\gamma$ we have $\gamma=\min_{x,y\in\R^+}S(x,y)$.
\end{theorem} 
In order to keep the proof of Thm.~\ref{thm:mainthm} reasonably concise, we will first start with two somewhat technical lemmas. In the first lemma, we will establish some periodicity properties of $S(x,y)$ and the group, in the case where we have multiple saddle points. 
\begin{lemma}\label{lemma:periodicity}
Let $(x_0,y_0)$ be a dominant saddle point, and $(x_i,y_i)=(\alpha_i x_0,\beta_ix_0)$ be the associated ones, with $S(x_i,y_i)=\zeta_iS(x_0,y_0)$. Define furthermore $\phi(x,y),\psi(x,y)$ such that $\Phi(x,y)=(x,\phi(x,y)),\Psi(x,y)=(\psi(x,y),y)$ are the generators of the group as in (\ref{eq:group1}),(\ref{eq:group2}). We then have, for all $x,y\in\mathbb{C}$:\begin{align}
\label{eq:periodicity1}
    S(\alpha_i x,\beta_i y)&=\zeta_i S(x,y),\\
    \label{eq:periodicity2}
    \psi(\alpha_i x,\beta_i y)&=\alpha_i \psi(x,y),\\
    \label{eq:periodicity3}
    \phi(\alpha_i x,\beta_i y)&=\beta_i \phi(x,y).
\end{align}
\end{lemma}
\textbf{Remark:} Lemma~\ref{lemma:periodicity} still holds true in more than two dimensions, with a completely analogous proof. Also note in particular that \eqref{eq:periodicity1}--\eqref{eq:periodicity3} hold true for all $x,y\in\C$, and not only saddle points.
\begin{proof}
    As $(x_0,y_0)$ is a dominant saddle point and $(x_i,y_i):=(\alpha_ix_0,\beta_i y_0)$ associated to it, we have $|\alpha_i|=|\beta_i|=|\zeta_i|=1$, and know that for each monomial $x^ky^l$ appearing in $S(x,y)$ we have $(\alpha_i x)^k(\beta_i y)^l=\zeta_i x^ky^l$. Consequently, for all such $k,l$ we have $\alpha_i^k\beta_i^l=\zeta_i$, and thus (\ref{eq:periodicity1}) holds.\\
    Since $\Psi$ by construction changes $(x,y)$ to another pair $(\psi(x,y),y)$ such that both $S(x,y)=S(\psi(x,y),y)$, we know (remember that $S(x,y)$ is quadratic in $x,y$) that if \begin{align}
        S(x_1,y)&= S(x_2,y),\\
        x_1&\neq x_2,
    \end{align}
   then we have $\psi(x_1,y)=x_2$. Now suppose we have a pair $(x,y)$ such that $\psi(x,y)\neq x$. Then, using (\ref{eq:periodicity1}), we can rewrite:\begin{align}\label{eq:periodicity4}
        S(\alpha_i \psi(x,y),\beta_i y)=\zeta_iS(\psi(x,y),y)=\zeta_i S(x,y)=S(\alpha_i x,\beta_i y).
    \end{align}
    In this case, from the above we can conclude that \eqref{eq:periodicity2} holds whenever $\psi(x,y)\neq y$. But $\psi(x,y)=x$ can be true only when $\frac{\partial S}{\partial x}(\cdot, y)=0$, i.e.~where for a given $y$ we have no solution other than $x$ to $S(\cdot,y)=S(x,y)$. This means that we have $\psi(x,y)\neq y$ almost everywhere. Since the functions on the left- and right-hand side of \eqref{eq:periodicity2} are rational and agree almost everywhere, we know that they are indeed the same. By a symmetric argument, we show \eqref{eq:periodicity3}. 
\end{proof}

Lastly, we will show that given an asymptotic representation as in Thm.~\ref{thm:mainthm} below, the functions appearing therein are indeed polyharmonic.
\begin{lemma}\label{lemma:mainlemma}
Suppose $q(B;n)$ is a (combinatorial) quantity satisfying \begin{align}
    \label{lemma1:countingcondition}
    q(B;n+1)=\sum_{s\in\mathcal{S}}\omega_s q(B-s;n)
\end{align}
for all $B\in\mathbb{Z}_{\geq 0}\times\mathbb{Z}_{\geq 0},n\geq 0$, and that at the same time it is of the form \begin{align}\label{eq:lemma_rep}
    q(B;n)=\frac{\gamma^n}{n^c} \left[\sum_{p=1}^{m-1} \frac{\sum_{i=1}^r v_{p,i}(B)\zeta_i^n}{n^{p}}+\mathcal{O}\left(\frac{1}{n^{m+1}}\right)\right],
\end{align}
for all $k\geq 0$, with the $\zeta_i$ pairwise different roots of unity. Then, the $v_{p,i}$ are $\gamma$-polyharmonic of degree $p$.\\
If, additionally, for a fixed point $B$ we know that $q(B;n)=0\text{ }\forall n\in\mathbb{N}$, then $v_{p,i}(B)=0$ for all $p,i$.

\end{lemma}

\begin{proof}
Substituting (\ref{eq:lemma_rep}) into (\ref{lemma1:countingcondition}) gives us \begin{align*}
    q(B;n+1)&=\sum_{s\in\mathcal{S}}\omega_s q(B-s;n)\quad\Leftrightarrow\\
    \gamma^{n+1} \left[\sum_{p=1}^{m-1} \frac{\sum_{i=1}^r v_{p,i}(B)\zeta_i^{n+1}}{(n+1)^{p+c}}+\mathcal{O}\left(\frac{1}{n^{m+c}}\right)\right]&=\gamma^n\sum_{s\in\mathcal{S}}\omega_s\left[\sum_{p=1}^{m-1}\frac{\sum_{i=1}^r v_{p,i}(B-s)\zeta_i^n}{n^{p+c}}\right].\end{align*}
    Extracting the terms for $p=1$ and noticing that the others are smaller by a factor of at least $\frac{1}{n}$, we obtain\begin{align}
    \label{eq:harmonicity:proof3}
    \gamma\frac{n^{c+1}}{(n+1)^{c+1}}\sum_{i=1}^rv_{1,i}(B)\zeta_i^{n+1}+\mathcal{O}\left(\frac{1}{n}\right)&=
    \sum_{s\in\mathcal{S}}\omega_s\sum_{i=1}^rv_{1,i}(B-s)\zeta_i^n.
\end{align}
Letting $n$ go to $\infty$, we have \begin{align}\label{eq:lemma1_1}
\gamma\sum_{i=1}^rv_{1,i}(B)\zeta_i^{n+1}+\mathcal{O}\left(\frac{1}{n}\right)=\sum_{s\in\mathcal{S}}\omega_s\sum_{i=1}^rv_{1,i}(B-s)\zeta_i^n.
\end{align}
All we need to do now is show that each of the $v_{1,i}(x)$ is $\gamma$-harmonic by itself, i.e.~that (\ref{eq:lemma1_1}) holds for each summation index separately. To do so, let us forget for a moment the part $\mathcal{O}\left(\frac{1}{n}\right)$ and solve the exact analogue of (\ref{eq:lemma1_1}). \newline
Since by assumption the $\zeta_i$ are all different, we know that the vectors \begin{align}
    (\zeta_1^m,\dots, \zeta_l^m),\quad 0\leq m\leq r-1
\end{align}
are linearly independent (written as a matrix, they give a Vandermonde matrix with determinant $\prod_{i< j}(\zeta_j-\zeta_i)\neq 0$). Therefore, the system \begin{align}\label{eq:harmonicity:proof2}
    \sum_{i=1}^r\left(c v_{i,1}(B)-\sum_{s\in \mathcal{S}}\omega_s v_{i,1}(B-s)\right)\zeta_i^n=0
\end{align}
has no nontrivial solutions, and hence \begin{align}\label{eq:harmonicity:proof1}
    c v_{i,1}(B)=\sum_{s\in\mathcal{S}}\omega_s v_{i,1}(B-s)\quad\forall 1\leq i\leq r.
\end{align}
All that remains to do now is to see that the error term $\mathcal{O}\left(\frac{1}{n}\right)$ in (\ref{eq:lemma1_1}) does not change this. To do so, suppose now that (\ref{eq:harmonicity:proof1}) is not satisfied for some $i$. Then we know that there are arbitrarily large $n$ such that (\ref{eq:harmonicity:proof2}) does not hold, i.e.~its right-hand side takes a value $\varepsilon_n\neq 0$. As the $\zeta_i$ are roots of unity, there are only finitely many values which $\varepsilon_n$ can take for different $n$, so we cannot have convergence of $\varepsilon_n$ to $0$. But then, choosing $n$ large enough, (\ref{eq:harmonicity:proof3}) cannot hold either; a contradiction. Thus, (\ref{eq:harmonicity:proof1}) must hold, and we know that the $v_{1,i}(B)$ are harmonic. By induction, applying the discrete Laplacian $\triangle$  to both sides of (\ref{eq:harmonicity:proof3}), we argue in the same fashion that the $v_{k,i}(B)$ must be polyharmonic of degree $k$.\newline
To show that $q(B;n)=0$ for all $n$ implies that $v_{p,i}(B)=0$ for all $n$, let us assume the opposite. So suppose that $q(B;n)=0$ for all $n$, but that at the same time we have $p,i$ such that $v_{p,i}(B)\neq 0$. Assume our $p$ to be minimal with this property. Then, we know that \begin{align}
    \sum_{i=1}^rv_{p,i}(B)\zeta_i^n=0,
\end{align}
because otherwise this would be a contradiction to (\ref{eq:lemma_rep}) for large $n$. But then we can utilize independence of the vectors $(\zeta_1^k,\zeta_2^k,\dots,\zeta_l^k)$ for $k=0,\dots,r-1$ as before and arrive at a contradiction.\\
Finally, to show that the $v_{p,i}$ are polyharmonic of order $p$, we can apply the operator $\triangle^{p-1}$ to both sides of (\ref{eq:lemma_rep}), notice that the first $p-1$ terms vanish by assumption and then repeat the above argument.\end{proof}
\textbf{Remarks:} 
\begin{itemize}
    \item In the proof of Thm.~\ref{thm:mainthm} we will see that in our case, for different $i$ the $v_{p,i}(k,l)$ differ only by a factor of $\alpha_i^{-k}\beta_i^{-l}$. Here, given a dominant saddle point $(x_0,y_0)$, and an associated one $(x_i,y_i)$, then $\alpha_i,\beta_i$ are the numbers such that $(x_i,y_i)=(\alpha_ix_0,\beta_i y_0)$. This allows us to essentially talk about only a single polyharmonic function $v_p$ for any given $p$; namely the one defined by the dominant saddle point.
\item By studying the exact shape of the functions $v_p$ more closely (in particular their degree as polynomials), we will see in Section~\ref{sec:startingpoints} that in the context of this article, we always have $\triangle v_{p+1}=v_p+r_{p-1}$, where $r_{p-1}$ is some polyharmonic function of degree at most $p-1$. This hints at the fact that the polyharmonic functions appearing in the asymptotic expansions are not arbitrary, but do have some form of recursive structure.
\end{itemize}

We now have all ingredients ready for the proof of Thm.~\ref{thm:mainthm}.

\begin{proof}[Proof (of Thm.~\ref{thm:mainthm})]
By the assumptions, we know that we have
\begin{align}
    q(0,(k,l);n)=\left[x^ky^lt^n\right]Q(x,y;t)=\left[x^{k+1}y^{l+1}t^n\right]\frac{\sum_{g\in\mathcal{G}}\operatorname{sgn}(g)g(xy)}{k(x,y;t)}.
    \end{align}
    As $k(x,y)=1-tS(x,y)$, we can rewrite \begin{align}
        q(0,(k,l);n)&=\left[x^{k+1}y^{l+1}t^n\right]\sum_{i=0}^\infty t^iS(x,y)^iN(x,y)\\
        &=\left[x^{k+1}y^{l+1}\right]S(x,y)^nN(x,y).
    \end{align}
    By Cauchy's formula, we have \begin{align}\label{eq:mainthm:proof1}
        q(0,(k,l);n)=-\frac{1}{4\pi^2}\int_{\Gamma_1}\int_{\Gamma_2}\frac{S(x,y)^nN(x,y)}{x^{k+2}y^{l+2}}\mathrm{d}x\mathrm{d}y,
    \end{align}
    with $\Gamma_{1,2}$ being closed curves around the origin. To evaluate the asymptotics of this integral, we utilize the saddle point method, as described for instance in~\cite[Chapter~VIII]{Flajolet}.\newline
    The main idea is to conveniently choose our contours $\Gamma_1,\Gamma_2$ such that they make the integral as easy to compute as possible.\\
    With this in mind, suppose that $(x_0,y_0)$ is a dominant saddle point, and pick $\Gamma_1=\{|x|=|x_0|\}, \Gamma_2=\{|y|=|y_0|\}$. We know that the modulus of $S(x,y)$ on $\Gamma_1\cap\Gamma_2$ is maximal; and the only other points where it attains the same value are the associated saddle points $(x_i,y_i)$. At any other point, $|S(x,y)|$ will be strictly smaller -- hence, when $n$ goes to infinity, it suffices to compute the integral locally around our saddle points, since the rest will grow exponentially slower. We could hypothetically run into issues if $N(x,y)$ were to be infinite, but we will see that this is not the case at our saddle points. For any other points, one can easily verify using the definition of the group that we can have infinite values of $N(x,y)$ only at lines of the form $\{x=\zeta_i\}$ or $\{y=\chi_i\}$ for finitely many $\zeta_i,\chi_i$. It follows from the residue theorem, applied for the two coordinates iteratively, that these singularities can safely be neglected. It is then easy to check that, given an $\varepsilon>0$, the set $\{(x,y):|S(x_0,y_0)-S(x,y)|<\varepsilon\}$ is contained in a domain of the form $\{(x,y):|x-x_0|<\varepsilon_1, |y-y_0|<\varepsilon_2\}$. As previously mentioned, in order to find the asymptotics, the rest of the integral can be neglected, as it will be exponentially smaller. Changing our coordinates to $x=x_0e^{is/\sqrt{n}},y=y_0e^{it/\sqrt{n}}$, this corresponds to a region of the form $\left|\frac{s}{\sqrt{n}}\right|<\delta_1, \left|\frac{t}{\sqrt{n}}\right|<\delta_2$, or, equivalently, $|s|<\delta_1\sqrt{n},|t|<\delta_2\sqrt{n}$.\newline
     To find the asymptotics, it therefore suffices to compute the integrals \begin{align}
        \int_{-\delta_1\sqrt{n}}^{\delta_1\sqrt{n}}\int_{-\delta_2\sqrt{n}}^{\delta_2\sqrt{n}}F_j(s,t,k,l,n)\mathrm{d}t\mathrm{d}s,
    \end{align}
    where $F_j(s,t,k,l,n)$ is the expression obtained by substituting $x=x_je^{is/\sqrt{n}},y=y_je^{it/\sqrt{n}}$, for $(x_j,y_j)$ the relevant saddle points (i.e.~$(x_0,y_0)$ and the ones associated to it, as outlined in Section~\ref{sec:prelims}). We will see that, given any fixed $m\in\mathbb{N}$, each such integral can be written in the form \begin{align}
        \frac{\gamma^n}{n^c}\int_{-\delta_1\sqrt{n}}^{\delta_1\sqrt{n}}\int_{-\delta_2\sqrt{n}}^{\delta_2\sqrt{n}}e^{-Q(s,t)}\left(p_0(s,t)+\frac{p_1(s,t)}{n}+\frac{p_2(s,t)}{n^2}+\dots+\mathcal{O}\left(\frac{1}{n^m}\right)\right)\mathrm{d}t\mathrm{d}s,
    \end{align}
    where $Q(s,t)$ is some (positive definite) quadratic form and the $p_j(s,t)$ are polynomials. Assume w.l.o.g. that $\delta_1\geq \delta_2$. For $s,t>\delta_2\sqrt{n}$, we can see that the integral over the remaining part of $\mathbb{R}$ is exponentially small in $n$, and therefore we can consider the integral \begin{align}\label{eq:mainthm:proof2}
        \frac{\gamma^n}{n^c}\int_{-\infty}^\infty\int_{-\infty}^\infty e^{-Q(s,t)}\left(p_0(s,t)+\frac{p_1(s,t)}{n}+\dots+\mathcal{O}\left(\frac{1}{n^m}\right)\right)\mathrm{d}t\mathrm{d}s.
    \end{align}
    Consequently, all we need to do is to consider the latter integral for all saddle points with maximum absolute value (of which there are finitely many), and by computing all the expressions up to a fixed $p_j(s,t)$, we will then have obtained the asymptotics of (\ref{eq:mainthm:proof1}) and at the same time shown (\ref{eq:mainthm:asymptotics}).\newline 
    In the following, we will proceed in two steps: first, we pick a dominant saddle point and show that, locally around this point, everything works out smoothly. Then, we pick any associated saddle point and show that, up to powers of roots of unity, nothing changes from the first step.
    \begin{enumerate}
    \item Suppose $(x_0,y_0)$ is a dominant saddle point, let $\gamma:=S(x_0,y_0)$ and fix an $m\in\mathbb{N}$. First, we show that $0=|N(x_0,y_0)|<\infty$. This is due to $(x_0,y_0)$ being a saddle point: we have $\frac{\partial S}{\partial y}(x_0,y_0)=0$ and thus $y_0$ is the unique solution to $S(x_0,\cdot)=S(x_0,y_0)$. Therefore, $\phi(y_0)=y_0$, and in the same manner we can see $\psi(x_0)=x_0$. It follows immediately that the alternating orbit sum of $xy$ evaluated at $(x_0,y_0)$ is $0$; in particular it is finite. \newline 
    Our next step is to show that we can rewrite the integrand (that is, the one in (\ref{eq:mainthm:proof1})) as in (\ref{eq:mainthm:proof2}). To do so, we substitute $x\mapsto x_0\exp{\frac{is}{\sqrt{n}}}=:e_s,y\mapsto y_0\exp{\frac{it}{\sqrt{n}}}=:e_t$, and then separate the integral into three parts: $S(e_s,e_t)^n, N(e_s,e_t)$ and the denominator $1/e_s^{k+1}e_t^{l+1}$ (note that a power in the denominator vanishes due to the substitution rule). 
    \begin{enumerate}
        \item \textbf{Part 1: $S(e_s,e_t)^n$}\\
    As $(x_0,y_0)$ is a dominant saddle point of $S(x,y)$ and therefore $\frac{\partial S}{\partial x}(x_0,y_0)=\frac{\partial S}{\partial y}(x_0,y_0)=0$, we have a Taylor expansion of $S(x,y)$ around $(x_0,y_0)$ of the form\begin{align}
    S(x,y)=\gamma-u(x-x_0)^2-v(x-x_0)(y-y_0)-w(y-y_0)^2+\dots
        \end{align}
    After our substitution, this gives us (note that $e_s,e_t$ are functions of $s,t$ and $n$) \begin{align}
        S\left(e_s,e_t\right)=\gamma-u\frac{s^2}{4n}-v\frac{st}{n}-w\frac{t^2}{4n}+A(s,t,n),
    \end{align}
    with $A(s,t,n)=n^{-3/2}\sum_{j=0}^{2m-1} \frac{a_j(s,t)}{n^{j/2}}+\mathcal{O}\left(n^{-m-3/2}\right)$, and the $a_j$ homogeneous polynomials. We know that $\widehat{Q}(s,t):=\left[\frac{u}{4}s^2+vst+\frac{w}{4}t^2\right]$ is a positive definite quadratic form (because our saddle point is a local minimizer of $S(x,y)$ in $\mathbb{R}^+\times\mathbb{R}^+$). Consequently, we can write \begin{align}
        \log S(e_s,e_t)&=\log\left(\gamma-\left[\widehat{Q}(s,t)-A(s,t,n)\right]\right)\\
        &=\log\gamma-\underbrace{\frac{1}{\gamma}\widehat{Q}(s,t)}_{=:Q(s,t)}+B(s,t,n),
    \end{align}
with once again $B(s,t,n)=n^{-3/2}\sum_{j=0}^{2m-1} \frac{b_j(s,t)}{n^{j/2}}+\mathcal{O}\left(n^{-m-3/2}\right)$, the $b_j$ homogeneous polynomials and $Q(s,t)$ a positive definite quadratic form. Consequently,\begin{align}
    S(e_s,e_t)^n=\exp\left[n\log S(e_s,e_t)\right]=\gamma^n\exp\left[-Q(s,t)\right]\exp\left[nB(s,t,n)\right].
\end{align}
Comparing this to (\ref{eq:mainthm:proof1}), the first two factors are already precisely as we want them, and the last factor is of the form $\sum_{m\geq 1}\frac{q_m(s,t)}{n^{m/2}}$, with the $q_m(s,t)$ homogeneous of degree $m$.
\item \textbf{Part 2: $N(e_s,e_t)$}\newline
As we have seen that $|N(x_0,y_0)|< \infty$, and as $N(x,y)$ is a rational function in $x,y$, it follows that we can write \begin{align}
    N(e_s,e_t)=\sum_{j=0}^{2m-1}\frac{d_j(s,t)}{n^{j/2}}+\mathcal{O}\left(\frac{1}{n^m}\right),
\end{align}
with the $d_j(s,t)$ homogeneous polynomials of degree $j$.
\item \textbf{Part 3: $1/e_s^{k+1}e_t^{l+1}$}\newline
Lastly, we have \begin{align}\label{eq:mainthm:degree}
    \frac{1}{e_s^{k+1}e_t^{l+1}}=x_0^{-k-1}y_0^{-l-1}\sum_{j\geq 0}\frac{1}{n^{j/2}}\frac{\left(-i[(k+1)s+(l+1)t]\right)^j}{j!}.
\end{align}
Note in particular that this last factor is the only part which depends on the endpoint $(k,l)$. 
\end{enumerate}
Multiplying the power series together and sorting them by powers of $n$, we obtain as a result the contribution of this saddle point to (\ref{eq:mainthm:proof2}) of the form \begin{align}\label{eq:mainthm:proof3}
    \frac{\gamma^n}{n^c}\left[\sum_{p=0}^{2m-1}\frac{1}{n^{p/2}}\frac{1}{x_0^{k+1}y_0^{l+1}}\int_{-\infty}^\infty\int_{-\infty}^\infty e^{-Q(s,t)}q'_{p}(s,t,k,l)\mathrm{d}s\mathrm{d}t+\mathcal{O}\left(\frac{1}{n^m}\right)\right],
\end{align}
where the $q'_{p}(s,t,k,l)$ are polynomials in $s,t,k,l$. One can easily see that $q'_p(s,t,k,l)$ is homogeneous of degree $p$; for odd $p$ the double integral therefore vanishes by symmetry. Thus we can rewrite (\ref{eq:mainthm:proof3}) as 
\begin{align}
        \frac{\gamma}{n^c}\left[\sum_{p=0}^{m-1}\frac{1}{n^{p}}\frac{1}{x_0^{k+1}y_0^{l+1}}\int_{-\infty}^\infty\int_{-\infty}^\infty e^{-Q(s,t)}q_{p}(s,t,k,l)\mathrm{d}s\mathrm{d}t+\mathcal{O}\left(\frac{1}{n^m}\right)\right],
\end{align}
with $q_p:=q'_{2p}$. The factor $\frac{1}{n^c}$ in (\ref{eq:mainthm:proof3}) stems from the fact that we obtain a factor of $1/n$ by the substitution rule, and that the integral might vanish for small values of $p$. In particular, from this we can conclude that the $c$ appearing in~\eqref{eq:mainthm:asymptotics} is integer.

\item Suppose now that we have another saddle point $(x_j,y_j)=(\alpha_j x_0,\beta_j y_0)$ associated to $(x_0,y_0)$, and pick $\zeta_j$ such that $S(x_j,y_j)=\zeta_jS(x_0,y_0)$ (notice that we then have $|\zeta_j|=1$, as discussed in Section~\ref{sec:prelims}). Our goal is now to describe the series expansion of the numerator around $(x_j,y_j)$ using the one around $(x_0,y_0)$. We substitute as before $x\mapsto e_s':=x_j\exp\frac{is}{\sqrt{n}},y\mapsto e_t':=y_j \exp\frac{it}{\sqrt{n}}$. Due to Lemma~\ref{lemma:periodicity}, we can now conclude that the series representation (w.r.t.~$n$ at $\infty$) around $S(e_s',e_t')$ is the same as the one of $\zeta_j S(e_s,e_t)$, and the representation of $N(e_s',e_t')$ is the same as $\alpha_j\beta_jN(e_s,e_t)$. Lastly, the expansion of $1/e_s^{\prime(k+1)}e_t^{\prime(l+1)}$ clearly changes only by adding a factor of ${\alpha_j}^{k+1}{\beta_j}^{l+1}$ as well. Therefore, we can conclude that the contribution of the saddle point $(x_j,y_j)$ is the same as the one of $(x_0,y_0)$ up to a factor of $\zeta_j^n\alpha_j^{-k}\beta_j^{-l}$.
\end{enumerate}
By Lemma~\ref{lemma:mainlemma}, we deduce that the $v_p(k,l)$ are indeed $\gamma$-polyharmonic of degree $p$ (note that for our dominant saddle point we have $\alpha=\beta=\zeta=1$).\newline
Because the only point in the construction where $k,l$ appear is in (\ref{eq:mainthm:degree}), one finds that the $v_p(k,l)$ are, up to a factor of $x_0^{-k}y_0^{-l}$, bivariate polynomials in $k,l$. Due to Lemma~\ref{lemma:saddlepoint_driftzero}, they are therefore polynomials if and only if the model has zero drift.\newline
The fact that $c=\pi/\theta$ (with $\theta$ the correlation coefficient as defined in (\ref{eq:deftheta})) follows from~\cite[Thm.~6]{Denisov}. To see that $\gamma=\min_{x,y\in\R^+}S(x,y)$, it suffices to remember that at the beginning of the proof we had $\gamma=S(x_0,y_0)$ for a dominant saddle point $(x_0,y_0)$, which is a minimizer of $S$ by its definition in Section~\ref{sec:prelims}.
\end{proof}

The construction used in the proof allows us to give some further properties of the polyharmonic functions appearing in the asymptotic expansion.

\begin{corollary}\label{cor:mainthm}
The degree of the polynomial part of $v_p(k,l)$ (that is, without the factor of $x_0^{-k}y_0^{-l}$) is $c+2p-1$. 
\end{corollary}
\begin{proof}
    By looking once again at the proof of Thm.~\ref{thm:mainthm} and in particular (\ref{eq:mainthm:degree}). From there, the statement follows immediately.
\end{proof}
\begin{corollary}\label{cor:mainthm_2}
For any orbit-summable model, we have $\pi/\theta\in\mathbb{N}$.
\begin{proof}
    From Thm.~\ref{thm:mainthm}, we know that $c=\pi/\theta$. From Cor.~\ref{cor:mainthm}, we know that $c+2p-1$ is the degree of the polynomial part of $v_p(k,l)$, and therefore integer. Since $2p-1$ is integer, the statement follows.
\end{proof}
\end{corollary}

\textbf{Remarks:} \begin{itemize}
    \item Thm.~\ref{thm:mainthm} holds true in higher dimensions as well, and indeed the proof translates directly. The one difference lies in the powers of $n$ which appear: by the substitution rule $\frac{\mathrm{d}x_i}{\mathrm{d}s_i}=c_ie^{is_i/\sqrt{n}}$, one obtains additional factors of $n^{-1/2}$. Thus, for even dimensions the constant $c$ in (\ref{eq:mainthm:asymptotics}) will be integer, whereas for odd dimensions it will be in $\frac{1}{2}+\mathbb{N}$. An example case for three dimensions is treated in~\ref{app:3dim}. 
    \item When looking at models with large steps, the one thing that could go wrong is that the numerator might have a singularity at a saddle point. Usually, this seems not to be the case, and for a given model this condition is very easy to check. An example is treated in~\ref{app:largesteps}.
    \item A table of the first three asymptotic terms for the $19$ unweighted orbit-summable models is given in~\cite[App.~D]{Nes23}.
    \end{itemize}
\subsection{Example: the Gouyou-Beauchamps model}\label{sec:GB_example}
In this section, we will illustrate the result of Thm.~\ref{thm:mainthm} by computing the asymptotics for the Gouyou-Beauchamps model, and in doing so find an explicit formula for the polyharmonic functions appearing therein.\newline 
The Gouyou-Beauchamps Walk is defined by the step set $\{\nwarrow,\searrow,\leftarrow,\rightarrow\}$. Its step polynomial is $S(x,y)=\frac{y}{x}+\frac{x}{y}+\frac{1}{x}+x$, and solving $\frac{\partial S}{\partial x}=\frac{\partial S}{\partial y}=0$ yields the four solutions $(x,y)=(\pm 1,\pm 1)$. We find that $S(1,1)=4, S(-1,1)=-4$, $S(1,-1)=S(-1,-1)=0$. Verifying that the second derivatives do not vanish, we therefore have the dominant saddle point $(1,1)$ and one other associated to it, namely $(-1,1)$, to consider. Note that the appearance of two saddle points is not at all surprising here, due to parity (or, in more general terms, periodicity) considerations: if the first coordinate of a given point is even, then we can only hit it after an even or odd number of steps. Therefore we can already expect at this point the asymptotics resulting from the saddle points to be precisely the same up to a factor of $(-1)^{k+n}$, which coincides with the statement of Thm.~\ref{thm:mainthm}. Hence, we will only consider the dominant saddle point here.\\
The alternating orbit sum of $xy$ can be checked to be \begin{align}\label{eq:exampleGB_N(x,y)}
    N(x,y):=-\frac{(-1 + x^2) (-1 + y) (x^4 + y^3 - x^2 y -x^2y^2)}{x^3 y^2}.
\end{align}
Out integrand is therefore of the form \begin{align*}
    \frac{S(x,y)^n N(x,y)}{x^{k+2}y^{l+2}},
\end{align*}
where after letting $x\mapsto e_s:=e^{is/\sqrt{n}}, y\mapsto e_t:=e^{it/\sqrt{n}}$, by the substitution rule we will end up with \begin{align*}
    -\frac{1}{n}\frac{S(e_s,e_t)^nN(e_s,e_t)}{e_s^{k+1}e_t^{l+1}}.
\end{align*}
The first factor of $-1/n$ is entirely harmless; we will therefore proceed to compute each of the factors in the second fraction separately, precisely as in the proof of Thm.~\ref{thm:mainthm}.

\begin{enumerate}\item\textbf{Series representation of $S(e_s,e_t)^n$:}\newline
First, we can compute (we let $m:=\sqrt{n}$ for readability) \begin{align}
        S(e_s,e_t)=4-\sum_{j\geq 2}\frac{a_j}{m^j}=4-\underbrace{\frac{s^2+(s-t)^2}{n}}_{=:4Q(s,t)}+\underbrace{\sum_{j\geq 3}\frac{a_j}{m^j}}_{=:A(s,t,m)},
    \end{align}
    with \begin{align*}
        a_j:=i^j\frac{s^j+(s-t)^j}{j!}(1+(-1)^j).
    \end{align*}
   From here we obtain \begin{align}
        \log\left[S(e_s,e_t)\right]&=\log 4-\frac{Q(s,t)}{n}+B(s,t,m),
    \end{align}
    where $B(s,t,m)=\sum_{j\geq 3}\frac{b_j}{m^j}$, with \begin{align*}
        b_j=\sum_{i_1+i_2+\dots+i_n=j}\frac{(-1/4)^n}{n}\prod_{p=1}^na_{i_p},
    \end{align*}
    with the $i_j$ positive integers.\newline
    Finally, we can compute \begin{align}
        S(e_s,e_t)^n&=\exp\left(n\log\left[S(e_s,e_t)\right]\right)\\
        &=\exp\left(n\log\left[4+B(s,t,m)\right]-Q(s,t)\right)\\
        &=4^ne^{-Q(s,t)}\underbrace{\exp\left[B(s,t,m)\right]}_{:=C(s,t,m)},
    \end{align}
    where $C(s,t,m)=\sum_{k\geq 1}\frac{c_j}{j^k}$ with \begin{align*}
        c_j=\sum_{i_1+i_2+\dots+i_n=j}\frac{1}{n!}\prod_{p=1}^nb_{i_p},
    \end{align*}
    with the $i_j$ again positive integers.

\item\textbf{Series representation of $N(e_s,e_t)$:}\newline
Using (\ref{eq:exampleGB_N(x,y)}), we find that \begin{align}\label{eq:Integrand_Part2}
        N\left(e_s,e_t\right)=\sum_{j\geq 1}\frac{d_j}{m^j},
    \end{align}
    with \begin{align*}
        d_j:=\frac{i^j}{j!}\left[-(-4)^js^j+(2s-3t)^j-2^j(s-t)^j-(t-2s)^j+(t-4s)^j+\right.\\\left.(-2)^j(s+t)^j-(-3)^jt^j\right].
    \end{align*}
    
    \item\textbf{Series representation of $1/e_s^{k+1}e_t^{l+1}$:}\newline
    Lastly, we have \begin{align*}
        \frac{1}{e_s^{k+1}e_t^{l+1}}&=\exp\left[-i(s(k+1)+t(l+1))/m\right]\\&=\sum_{j\geq 0}\frac{1}{m^j}\underbrace{\frac{\left(-i((k+1)s+(l+1)t)\right)^j}{j!}}_{=:f_j}.
    \end{align*}
\end{enumerate}
    Overall, we obtain as product of the three factors computed above \begin{align*}
        4^n e^{-Q(s,t)} \sum_{p\geq 0}\frac{1}{m^p}\underbrace{\left[\sum_{j_1+j_2+j_3=p}c_{j_1}d_{j_2}f_{j_3}\right]}_{=:q_p},
    \end{align*}
    for $j_1,j_2,j_3$ nonnegative integers. In particular, we notice that $q_p$ is homogeneous of degree $p$ in $s,t$, and of degree $p$ in $k,l$. In order to compute the contribution up to order $\mathcal{O}\left(\frac{1}{n^j}\right)$ of this saddle point to the asymptotics, all we need to do now is compute \begin{align*}
        4^n\frac{1}{4\pi^2n}\sum_{p=0}^{2r-1}\int_{-\infty}^\infty\int_{-\infty}^\infty e^{-Q(s,t)}q_p(s,t,k,l)\mathrm{d}t\mathrm{d}s. 
    \end{align*}
    This gives us an explicit formula for the asymptotics of this model. In particular, we can check directly that all coefficients of $\frac{1}{m^k}$ for odd $k$ vanish. By computing the integrals, we see that \begin{small}\begin{align*}
        v_1(k,l)=&\frac{64}{\pi} (1 + k) (1 + l) (2 + k + l) (3 + k + 2 l),\\
        v_2(k,l)=&-\frac{32}{\pi} (1 + k) (1 + l) (2 + k + l) (3 + k + 2 l) (35 + 2 k^2 + 
   4 k (2 + l) + 4 l (3 + l)),\\
   v_3(k,l)=&\frac{8}{\pi} (1 + k) (1 + l) (2 + k + l) (3 + k + 2 l) (25 + 2 k^2 + 4 k (2 + l) \\&+ 4 l (3 + l)) (61 + 2 k^2 + 4 k (2 + l) + 
   4 l (3 + l)).
    \end{align*}
    \end{small}
\noindent The harmonic function $v_1(k,l)$ was already computed some time ago in~\cite{Bia92,Ras11,CMMR17}.

\subsection{Periodicity}\label{sec:periodicity}
Considering the combinatorial context, it is clear that in any asymptotic expansion as in (\ref{eq:mainthm:asymptotics}), the discrete harmonic function $v_1(k,l)$ will always be positive (if it does not vanish). When looking at the computations in \cite[App.~D]{Nes23}, it appears as if there is an even stronger pattern; namely that $v_p(k,l)$ is positive for even $p$, and negative for odd $p$. It turns out, however, that this is not generally true; a counterexample is given for instance by computing enough terms in the expansion of the Simple Walk~\cite{Chapon_SW}.\\
It is also a direct consequence of Thm.~\ref{thm:mainthm} that the number of saddle points is closely tied to the periodicity of the model. If we have a single saddle point, then clearly our model is aperiodic; but the number of saddle points also corresponds directly to the periodicity.
\begin{lemma}
\label{lemma:sp_periodicity}
    Suppose that our model is irreducible (that is, the semi-group generated by step set $\mathcal{S}$ is all of $\mathbb{Z}^2$), and that it is $m$-periodic. Then we have exactly $m-1$ saddle points $s_1,\dots, s_{m-1}$ associated to our dominant saddle point $s_0=(x_0,y_0)$. The $\zeta_i$ corresponding to these saddle points (see Section~\ref{sec:prelims}) are -- in some order -- the $m$-th roots of unity.
\end{lemma}
\begin{proof}
    Let $r$ be the number of saddle points. We use the representation (\ref{eq:mainthm:asymptotics}) for the asymptotics of $q(k,l;n)$. We know that the $\alpha_i,\beta_i,\zeta_i$ are roots of unity. We can therefore pick $k,l$ such that $\alpha_i^k=\beta_i^l=1$ for all $i$. We then know, since our model is $m$-periodic, that there is a $z$, $0\leq z\leq m-1$, such that we have \begin{align}
        \chi(n):=1+\zeta_1^n+\dots+\zeta_{r-1}^n=\begin{cases}
            0\quad \text{  }n\not\equiv z\mod{m},\\
            m\quad n\equiv z\mod{m}.
        \end{cases}
    \end{align}
    Note that as value for $\chi(k\cdot m),k\in\mathbb{Z}$, we could pick any constant, because to compensate we can just multiply the corresponding polyharmonic functions in (\ref{eq:mainthm:asymptotics}) with a constant factor. By definition, we know that $\chi$ is a character on $\mathbb{Z}/m\mathbb{Z}$. Therefore, it can be uniquely written as sum of irreducible characters~\cite{Mas98}, which in this case are all the $m$-th roots of unity. From this, and the fact that the $\zeta_i$ are pairwise different, it follows already that the $\zeta_i$ are (in some ordering) the $m$-th roots of unity. 
\end{proof}
The irreducibility condition in Lemma~\ref{lemma:sp_periodicity} is necessary, as can be seen for the Diagonal Walk for instance. This walk is $2$-periodic, but we still have $4$ saddle points. The reason for the two extra saddle points is that there are some points the walk will never reach, which, heuristically speaking, translates to two additional conditions on $k,l$ for which two saddle points then do not suffice to express them.

\section{Quadrant walks with arbitrary starting point}
\label{sec:quadrantwalks_general}
\subsection{Full asymptotic expansion}
\label{sec:startingpoints}
If we want to count walks starting from an arbitrary point $(k,l)$, then the only thing that changes is that we have a different monomial $x^{k+1}y^{l+1}$ in the functional equation (\ref{eq:prelims:feq}). It follows that we can proceed in exactly the same manner as before in order to obtain an expression as in (\ref{eq:prelims:osummable}), where only the sum on the right hand side changes. Hence, we can proceed in the same manner as for Thm.~\ref{thm:mainthm}, which allows us to recover the result of~\cite[Thm.~6]{Denisov}, which states that the first order term in the asymptotics of the number of walks terminating at a point $B=(u,v)$ and starting at a point $A=(k,l)$ will be given -- up to, once again, an exponential term and some power of $n$, and possible parity constraints -- by the product of two functions; a harmonic function in $B$ and a function in $A$ which is adjoint harmonic. This is fairly natural because the underlying combinatorial problem is highly symmetrical: 
any path from $A$ to $B$ corresponds to a path from $B$ to $A$ with reversed steps. This, and a similar statement for the higher order terms, can be formalized in Thm.~\ref{thm:startingpoint} below and the following Thm.~\ref{thm:startingpoint2}.
\begin{theorem}
    \label{thm:startingpoint}
    Suppose that $\mathcal{S}$ is a step set satisfying the general assumptions stated in Section~\ref{sec:prelims} and that $\mathcal{S}$ is orbit-summable. Then, with $(x_0,y_0)$ being a dominant saddle point, then there are constants $c,r\in\mathbb{N}$, $\gamma\in\R^+$, and $\alpha_i,\beta_i,\zeta_i\in\C$ as well as functions $v_p(k,l,u,v)$ such that for any $m\in\mathbb{N}$ we have \begin{align}
    q((u,v),(k,l);n)=\frac{\gamma^n}{n^c}\left[\sum_{p=1}^{m-1}\frac{v_p(k,l,u,v)\sum_{i=1}^r\alpha_i^{u-k}\beta_i^{v-l}\zeta_i^n}{n^p}+\mathcal{O}\left(\frac{1}{n^m}\right)\right].
    \end{align}
    The $v_p(k,l,u,v)$ are polynomials precisely if the drift is zero (else they contain exponential factors). In this case they are of bidegree $c+2p-1$ in both $(k,l)$ and $(u,v)$, and of total degree $2c+4p-2$. Each $v_p(k,l,u,v)$ is multivariate polyharmonic of degree $p$.
    
\end{theorem}

\textbf{Remark:} the constants $\alpha_i,\beta_i,\zeta_i$ stem from the saddle points associated to the dominant one via $s_i=(x_i,y_i)=(\alpha_ix_0,\beta_iy_0)$. Also note that despite the appearance of complex numbers, all sums end up being real.\newline 

Before proving Thm.~\ref{thm:startingpoint}, we first show the following lemma, which is a natural extension of Lemma~\ref{lemma:mainlemma}.

\begin{lemma}
    \label{lemma:multiharmonic}
    Suppose that $q(A,B;n)$ is a (combinatorial) quantity satisfying \begin{align}
    \label{eq:lemma_multharm_1}
        \sum_{s\in\mathcal{S}}\omega_sq(A-s,B;n-1)=q(A,B;n),\\
        \label{eq:lemma_multharm_2}
        \sum_{s\in\mathcal{S}}\omega_sq(A,B+s;n-1)=q(A,B;n),
    \end{align}
    for all $A,B\in\Z_{\geq}\times\Z_{\geq0}$, $n\geq 0$. Assume furthermore that $q(A,B;n)$ has an asymptotic representation of the form \begin{align}
        q(A,B;n)=\frac{\gamma^n}{n^c}\left[\sum_{p=1}^{m-1}\frac{\sum_{i=1}^lv_{p,i}(A,B)\zeta_i^n}{n^{p}}+\mathcal{O}\left(\frac{1}{n^{m}}\right)\right]
    \end{align}
    for all $m$, with the $\zeta_i$ pairwise different and of modulus $1$. Then, each $v_{p,i}$ is multivariate polyharmonic of order $p$.
\end{lemma}

\begin{proof}
The proof works in the very same manner as the proof of Lemma~\ref{lemma:mainlemma}; at each step we can choose whether to apply the identity (\ref{eq:lemma_multharm_1}) or (\ref{eq:lemma_multharm_2}), leading to an additional instance of $\triangle$ or $\tilde{\triangle}$ respectively.
\end{proof}

\begin{proof}[Proof (of Thm.~\ref{thm:startingpoint})]
    Analogous to the proof of Thm.~\ref{thm:mainthm}. We note that the contribution of $N(e'_s,e'_y)$ for an associated saddle point changes by a factor of $\alpha_i^{u+1}\beta_i^{v+1}$ instead of a factor of $\alpha_i\beta_i$, and only need to keep track of the coefficients $u,v$ throughout, which however behave exactly in the same manner as $k,l$.\\
    Finally, the polyharmonicity properties are a direct consequence of Lemma~\ref{lemma:multiharmonic}.
\end{proof}

\textbf{Remark:}
\begin{itemize}
    \item While the starting point $(u,v)$ and the end point $(k,l)$ of the walk end up playing a very similar role (see also Thm.~\ref{thm:startingpoint2}), which is not at all surprising from a combinatorial point of view, this is not at all obvious from the proof: the role of $(k,l)$ is very easily summarized as these coefficients only appear in the integrand as a factor of $x^{-k-1}y^{-l-1}$, the starting point does not appear directly as factor $x^{u+1}y^{v+1}$, but instead as its orbit sum. A priori, without the combinatorial interpretation, it does not seem to be obvious that both of these occurrences lead to a symmetrical role in the result.
\end{itemize}

In the following, we will want to describe the coefficients $v_p(k,l,u,v)$ appearing in Thm.~\ref{thm:startingpoint} more precisely. We will consider the drift zero case, which is however not a real restriction as we can transform any other model to a zero drift one using the Cramer transformation discussed in Section~\ref{sec:prelims}. The goal will be to prove the following theorem:\begin{theorem}
    \label{thm:startingpoint2}
    The polynomials $v_p(k,l,u,v)$ in the zero drift case of Thm.~\ref{thm:startingpoint} each have a representation of the form \begin{align}\label{eq:thm:startingpoint2}
        v_p(k,l,u,v)=
        \sum_{\substack{1\leq i,j\leq p, \\ i+j\leq p+1}}
        a_{i,j}h_i^j(k,l)g_i^j(u,v),
    \end{align}
   where the $a_{i,j}$ are constants, the $h_i^j$ are polyharmonic of degree (at most) $i$, and the $g_i^j$ are adjoint polyharmonic of degree (at most) $p+1-i$.
\end{theorem}

In order to prove Thm.~\ref{thm:startingpoint2}, it turns out to be very useful to have a polynomial basis of the space of polyharmonic functions for any given model. As the group is finite by assumption and $\pi/\theta\in\Z$ by Cor.~\ref{cor:mainthm_2}, we can use the basis given in~\cite{Nes22}, which consists of sequences $h_n^m$ of $n$-polyharmonic functions satisfying:\begin{enumerate}
    \item $\triangle h_{n+1}^m=h_n^m$,
    \item the $h_n^m(k,l)$ are bivariate polynomials of increasing degree in both $n$ and $m$: the degree will increase by $2$ for each step in $n$, whereas it will increase by at least $2$ and at most $c+1$ for each step in $m$. 
\end{enumerate}

\begin{proof}[Proof (of Thm.~\ref{thm:startingpoint2})]
    Taking $u,v$ as parameters, we can for each $(u,v)$ write $v_p(k,l,u,v)$ as a sum of the basis functions $h_n^m(k,l)$, and obtain \begin{align}
        \sum_{\substack{1\leq n,m\leq p, \\ n+m\leq p+1}}h_n^m(k,l)g_n^m(u,v).
    \end{align}
    Since we know that $v_p(k,l,u,v)$ is a bivariate polynomial of bidegree $c+2p-1$ in both $(k,l)$ and $(u,v)$ (which is also where the conditions $n,m\leq p$ and $n+m\leq p+1$ come from), the only thing we need to show is that $g_n^m$ is adjoint polyharmonic of degree $p+1-n$ for any $m$. To do so, we utilize Lemma~\ref{lemma:multiharmonic}. First, consider $n=p$. We then have \begin{align}
        \tilde{\triangle}\left(\triangle^{p-1}q_p(k,l,u,v)\right)=\tilde{\triangle}h_1^1(k,l)g_p^1(u,v)=h_1^1(k,l)\tilde{\triangle}g_p^1(u,v)=0,
    \end{align}
    therefore $v_p^1(u,v)$ is adjoint harmonic. As the discrete Laplacians are linear, we can now proceed by induction, in each step applying the same argument as above to all terms which are polyharmonic (in $(k,l)$) of order $n$, i.e.~the multiples of $h_n^1,\dots, h_n^{p+1-n}$. The statement follows. 
    \end{proof}

\textbf{Remarks:}
\begin{itemize}
\item Thm.~\ref{thm:startingpoint2} tells us in particular that we can write each coefficient $v_p$ as a sum of products of polyharmonic and adjoint polyharmonic functions, so that the degree of the former and latter adds up to at most $p+1$. This can be viewed as an extension of~\cite[Thm.~6]{Denisov}, where it is shown that $v_1$ is the product of a harmonic and an adjoint harmonic function (albeit in a much more general setting). 
\item By simple degree considerations, the only base function $h_p^j$ appearing in $v_p(k,l,u,v)$ can be $h_p^1$. From this it follows that $\triangle v_{p+1}(k,l,u,v)=v_p(k,l,u,v)+r_{p-1}(k,l,u,v)$, with $r_{p-1}$ polyharmonic of degree at most $p-1$. 
\item If the step set is symmetric in the sense that $\omega_{s}=\omega_{-s}$, then one can easily see that $a_{i,j}=a_{j,i}$ due to the symmetry of the underlying combinatorial problem ($q(A,B;n)=\tilde{q}(B,A;n)$, where $\tilde{q}$ denotes the paths with reversed steps). This holds true for an appropriate choice of basis in some other cases as well; for examples of this, see~\cite[App.~C]{Nes23}.
\item A priori it is not at all clear which elements of the basis $h_i^j$ constructed in~\cite{Nes22} actually appear in connection with some combinatorial problems, i.e.~if this basis is combinatorially reasonably chosen. As we will see in Section~\ref{sec:startingpoints:SW}, this seems to be the case.
\item We can in fact give an upper bound on the number of summands appearing in the decomposition (\ref{eq:thm:startingpoint2}). Denote by $\left(h_n^m\right)$ the basis of harmonic functions from \cite{Nes22} as discussed above. First of all, we can see by degree considerations that for any given $p$, all $(n,m)$ such that $h_n^m$ can appear in the decomposition are contained in the subset $\{(n,m):n+m\leq p+1\}$. If one then writes these $h_n^m$ as a table, this gives us a triangular shape of size increasing with $p$. Along the lines $n+m=k$, we will have functions of degree at least $2(k-1)$. We can now count the number of possible products of a given function $h_i^j$ with a base function $\tilde{h}_m^n$ of the adjoint Laplacian. In order not to exceed the degree of $v_p$, we can multiply $h_1^1$ with the entire triangle of adjoint base functions. For $h_1^2$ and $h_2^1$, we cannot multiply them with any $\tilde{h}_m^n$ with $m+n=p+1$, and so on. All in all, this gives us a maximum of $\frac{p(p+1)(p+2)(p+3)}{24}$ summands. This maximum is achieved e.g.~for the Simple Walk for $v_{1,2,3}$. Generally, the larger the value of $\pi/\theta$, the less summands we will have (since the degrees of the $h_i^j$ increase more quickly). 
\end{itemize}

In the following, we will see what this decomposition looks like in case of the Simple Walk. The Gouyou-Beauchamps model and the Tandem Walk are treated in~\cite[App.~C]{Nes23}.
\subsection{Example: the Simple Walk}\label{sec:startingpoints:SW}
In the following, to keep the expressions a bit shorter, we will give the expressions as after the substitution $k\mapsto k-1,l\mapsto l-1$ etc. (i.e. we have $kl$ instead of $(k+1)(l+1)$). This corresponds to a shift of the quarter plane, where instead of $\mathbb{Z}_{\geq 0}\times\mathbb{Z}_{\geq 0}$ we now consider $\mathbb{Z}_{>0}\times\mathbb{Z}_{>0}$. For the Simple Walk, with $\mathcal{S}=\{\rightarrow,\downarrow,\leftarrow,\uparrow\}$, we then have (after rescaling by multiplicative constants)\begin{align*}
h_1^1(k,l)=&kl,\\
h_1^2(k,l)=&k l (k - l) (k + l),\\
h_1^3(k,l)=&k l (14 - 5 k^2 + 3 k^4 - 5 l^2 - 10 k^2 l^2 + 3 l^4),\\
   h_2^1(k,l)=&kl(l-1)(l+1),\\
   h_2^2(k,l)=&kl (l-1) (l+1) (7 + 5 k^2 - 3 l^2),\\
   h_3^1(k,l)=&k l(l-2) (l-1)(l+1) (l+2).
\end{align*}
By symmetry, we can pick the base functions $\tilde{h}_i^j(u,v)=h_i^j(u,v)$ for the adjoint Laplacian. For the first three asymptotic terms with arbitrary starting and ending points, we obtain (again up to multiplicative constants)\begin{align*}
    v_1(k,l,u,v)=&kluv,\\
    v_2(k,l,u,v)=& k l u v (7 + 2 k^2 + 2 l^2 + 2 u^2 + 2 v^2),\\
    v_3(k,l,u,v)=&k l u v (167 + 140 k^2 + 12 k^4 + 140 l^2 + 24 k^2 l^2 + 
   12 l^4 + 140 u^2 + 40 k^2 u^2\\& + 24 l^2 u^2 + 12 u^4 + 140 v^2 + 
   24 k^2 v^2 + 40 l^2 v^2 + 24 u^2 v^2 + 12 v^4).
\end{align*}

One can check that Cor.~\ref{thm:startingpoint2} takes the form \begin{align*}
    v_1=&h_1^1\tilde{h_1^1},\\
    v_2=&4\left(h_2^1\tilde{h_1^1}+h_1^1\tilde{h_2^1}\right)+2\left(h_1^2\tilde{h_1^1}+\tilde{h_2^1}h_1^1\right)+15h_1^1\tilde{h_1^1},\\
    v_3=&\frac{192}{5}\left(h_3^1\tilde{h_1^1}+\tilde{h_3^1}h_1^1\right)+\frac{64}{5}\left(h_2^2\tilde{h_1^1}+\tilde{h_2^2}h_1^1\right)+4\left(h_1^3\tilde{h_1^1}+\tilde{h_1^3}h_1^1\right)+64\left(h_2^1\tilde{h_1^2}+\tilde{h_2^1}h_1^2\right)\\&+128h_2^1\tilde{h_2^1}+24h_1^2\tilde{h_1^2}+576\left(h_2^1\tilde{h_1^1}+\tilde{h_2^1}h_1^1\right)+288\left(h_1^2\tilde{h_1^1}+\tilde{h_1^2}h_1^1\right)+951h_1^1\tilde{h_1^1}.
\end{align*}
As the degree of $h_i^j$ is truly increasing by only $2$ whenever we increase either $i$ or $j$ by one, it turns out that we have indeed $1,5$ and $15$ different summands respectively. Due to the symmetry of this model, the second and third equations can be simplified a bit: for $v_2$, letting $g_2:=4h_2^1+2h_1^2$ (clearly, $g_2$ is then biharmonic) gives us \begin{align*}
    v_2=g_2\tilde{h_1^1}+h_1^1\tilde{g_2}+15h_1^1\tilde{h_1^1}.
\end{align*}
For $v_3$, letting in the same manner $g_3:=\frac{192}{5}h_3^1+\frac{64}{5}h_2^2+4h_1^3$, we have \begin{align*}
    g_3\tilde{h_1^1}+\tilde{g_3}h_1^1+144 \left(g_2\tilde{h_1^1}+\tilde{g_2}h_1^1\right)+64\left(h_2^1\tilde{h_1^2}+\tilde{h_2^1}h_1^2\right)+951h_1^1\tilde{h_1^1}.
\end{align*}
While one can view the definition of $g_{2}$ and $g_3$ as purely a crutch to make the resulting expressions shorter, they do in fact give in a sense a natural decomposition of the $v_i$: $g_2$ consists of the highest order terms in $v_2$, while $g_3$ consists of those in $v_3$.\hfill

\vspace{3mm}
\textbf{Remarks:} \begin{itemize}
    \item By the above, when taking a scaling limit, then we have \begin{align*}\lim_{\mu\to 0} \mu^{\alpha_i} v_i\left(\frac{x}{\mu},\frac{y}{\mu}\right)=\lim_{\mu\to 0}\mu^{\alpha_i} g_i\left(\frac{x}{\mu},\frac{y}{\mu}\right),
    \end{align*}
where $\alpha_i$ is an appropriate scaling constant, which means that the $g_i$ already give us all the terms which will not vanish in this kind of limit.
\item 
By~\cite[Thm.~2.3]{Poly}, we know that the continuous heat kernel $p_t(x,y,u,v)$ of a Brownian motion with covariance matrix \begin{align*}
    \Sigma=\begin{pmatrix}\sigma_{11}& \sigma_{12}\\\sigma_{12} & \sigma_{22}\end{pmatrix}
\end{align*}
with $\sigma_{11}=\mathbb{E}[X^2], \sigma_{12}=\mathbb{E}[XY], \sigma_{22}=\mathbb{E}[Y^2]$ (the scaling limit of this model) allows for an asymptotic representation of the form \begin{align}
p_t(x,y,u,v)=\frac{1}{t^2}\sum_{k\geq 1}\frac{f_k(x,y,u,v)}{t^k}.
\end{align}
Noticing that this representation looks almost the same as the one for the discrete case in Thm.~\ref{thm:startingpoint}, it is natural to compare the functions $v_p$ to their continuous counterparts $f_p$. For $v_1$, for instance, by~\cite[Lemma~13]{Denisov} we know that we will have (after appropriate scaling by a constant) $v_1\to f_1$. However, this is not at all clear for $p>1$. For the Simple Walk, one can check that $v_2\to f_2$, but this fails for $p=3$: we have \begin{align*}f_3(k,l,u,v)=k l u v (3 k^4 + 6 k^2 l^2 + 3 l^4 + 22 k^2 u^2 \\+ 6 l^2 u^2 + 3 u^4 + 
   6 k^2 v^2 + 22 l^2 v^2 + 6 u^2 v^2 + 3 v^4),\end{align*}
   whereas the scaling limit of $v_3(k,l,u,v)$ turns out to be
 \begin{align*}
       k l u v (3 k^4 + 6 k^2 l^2 + 3 l^4 + 10 k^2 u^2 \\+ 6 l^2 u^2 + 3 u^4 + 
   6 k^2 v^2 + 10 l^2 v^2 + 6 u^2 v^2 + 3 v^4),
   \end{align*}

   where the coefficients of $k^2u^2$ and $l^2v^2$ do not match. While it might seem a bit disheartening that the asymptotics of the discrete case do not converge to those of the continuous case `arbitrarily well', one could just as well argue that this is in fact not something that can be expected: \cite{Denisov} only gives us first-order convergence, so in a sense anything going beyond the first-order terms is already more than we could've hoped for. 
\end{itemize}

\section{Outlook}
\begin{itemize}
\item We know by~\cite[Thm.~6]{Denisov} that, for any (not necessary orbit-summable) model with a finite number of steps, the first order term of the asymptotics of $q(k,l;n)$ behaves as in Thm.~\ref{thm:mainthm}. However, not much is known about the higher order terms. In this context, one could view Thm.~\ref{thm:mainthm} as a partial generalization for the case of orbit-summable models with finite group, but it would be interesting to see whether something similar holds true for a more general class of models, which could for example also include models with large steps, as in e.g.~\cite{FR15}. While it appears at first that this might fail due to the more complicated structure of the group, in~\cite{BBMM21} the notion of an orbit is extended to a more general context, which might make an approach using a saddle point method feasible.
\item While in the zero-drift case there is a unique combinatorially relevant harmonic function for a given model (namely the positive one), there is no equivalent in the polyharmonic case, nor a reasonable combinatorial interpretation. It would be very interesting to know if there is such an interpretation, or at least, in some sense, a 'canonical' $p$-polyharmonic function. Intuitively, one may think such a canonical function should coincide in its scaling limit with the (continuous) polyharmonic function appearing in the corresponding continuous heat kernel, but as the higher order asymptotics of the discrete and continuous cases do not coincide (see Section~\ref{sec:startingpoints:SW}), this might not be an ideal choice. Maybe the functions $g_2,g_3,\dots$ as defined at the end of Section~\ref{sec:startingpoints:SW} could serve as candidates for such canonical representatives instead. 
\item While the saddle point method as used in this article works well to compute precise asymptotics for any given starting point given that the model has finite group and is orbit-summable, it seems that it would be hard to apply when one of these conditions does not hold. If the path counting function is algebraic (which, at least in the unweighted case, implies a finite group), sometimes it is possible to obtain an explicit expression which can then be utilized to extract asymptotics \cite{BM05,BM16_Gessel,FH84}. In the infinite group case, however, things seem to be more complicated. In the one-dimensional setting, a probabilistic approach seems to work in order to show an asymptotic expansion similar to (\ref{eq:intro:asymptoticseries}), using only moment conditions~\cite{DTW23}. \newline In two dimensions, one can tackle this problem using an explicit parametrization of the zero-set of the kernel via elliptic functions (in particular Jacobi theta-functions as in~\cite{EP20}). This method can then be utilized to show for some cases with infinite group (for instance the model with steps $\{\leftarrow,\uparrow,\rightarrow,\downarrow,\nearrow\}$) that the asymptotics of the number of walks returning to the origin contain logarithmic terms~\cite{KilianAndrewPrivate}, which makes an asymptotic expansion as in (\ref{eq:mainthm:asymptotics}) impossible. However, it turns out that the coefficients still have a similar structure in the sense that they consist of polyharmonic functions. 

\item Finally, one could pose similar questions for different domains, be it a higher dimension or a different cone, for instance the three quarter plane. See for example~\cite{BM16,Tro22,RT19,BBMM21, BM21, Yat17,BBKM16}.

\end{itemize}

 \bibliographystyle{elsarticle-harv} 
 %%alternatives: \bibliographystyle{elsarticle-num}
 %% \bibliographystyle{elsarticle-num-names}
 %%both not alphabetically ordered...
 \bibliography{cas-refs}

\begin{thebibliography}{50}
\expandafter\ifx\csname natexlab\endcsname\relax\def\natexlab#1{#1}\fi
\providecommand{\url}[1]{\texttt{#1}}
\providecommand{\href}[2]{#2}
\providecommand{\path}[1]{#1}
\providecommand{\DOIprefix}{doi:}
\providecommand{\ArXivprefix}{arXiv:}
\providecommand{\URLprefix}{URL: }
\providecommand{\Pubmedprefix}{pmid:}
\providecommand{\doi}[1]{\href{http://dx.doi.org/#1}{\path{#1}}}
\providecommand{\Pubmed}[1]{\href{pmid:#1}{\path{#1}}}
\providecommand{\bibinfo}[2]{#2}
\ifx\xfnm\relax \def\xfnm[#1]{\unskip,\space#1}\fi
%Type = Inbook
\bibitem[{Asmar and Grafakos(2018)}]{HFApplications}
\bibinfo{author}{Asmar, N.H.}, \bibinfo{author}{Grafakos, L.}, \bibinfo{year}{2018}.
\newblock \bibinfo{title}{Harmonic Functions and Applications}. \bibinfo{publisher}{Springer International Publishing}, \bibinfo{address}{Cham}. chapter~\bibinfo{chapter}{6}.
\newblock pp. \bibinfo{pages}{367--402}.
\newblock \DOIprefix\doi{10.1007/978-3-319-94063-2\_6}.
%Type = Article
\bibitem[{Banderier and Flajolet(2002)}]{BF02}
\bibinfo{author}{Banderier, C.}, \bibinfo{author}{Flajolet, P.}, \bibinfo{year}{2002}.
\newblock \bibinfo{title}{Basic analytic combinatorics of directed lattice paths}.
\newblock \bibinfo{journal}{Theor. Comput. Sci.} \bibinfo{volume}{281}, \bibinfo{pages}{37--80}.
\newblock \DOIprefix\doi{10.1016/S0304-3975(02)00007-5}.
%Type = Incollection
\bibitem[{Banderier and Wallner(2019)}]{BW19}
\bibinfo{author}{Banderier, C.}, \bibinfo{author}{Wallner, M.}, \bibinfo{year}{2019}.
\newblock \bibinfo{title}{The kernel method for lattice paths below a line of rational slope}, in: \bibinfo{booktitle}{Lattice path combinatorics and applications. Based on the 8th international conference on lattice path combinatorics and applications, California State Polytechnic University, Pomona (Cal Poly Pomona), CA, USA, August 17--20, 2015}. \bibinfo{publisher}{Cham: Springer}, pp. \bibinfo{pages}{119--154}.
\newblock \DOIprefix\doi{10.1007/978-3-030-11102-1\_7}.
%Type = Article
\bibitem[{Ba{\~n}uelos and Smits(1997)}]{BS97}
\bibinfo{author}{Ba{\~n}uelos, R.}, \bibinfo{author}{Smits, R.G.}, \bibinfo{year}{1997}.
\newblock \bibinfo{title}{Brownian motion in cones}.
\newblock \bibinfo{journal}{Probab. Theory Relat. Fields} \bibinfo{volume}{108}, \bibinfo{pages}{299--319}.
\newblock \DOIprefix\doi{10.1007/s004400050111}.
%Type = Incollection
\bibitem[{Biane(1992)}]{Bia92}
\bibinfo{author}{Biane, P.}, \bibinfo{year}{1992}.
\newblock \bibinfo{title}{Minuscule weights and random walks on lattices}, in: \bibinfo{booktitle}{Quantum probability and related topics}. \bibinfo{publisher}{Singapore: World Scientific}, pp. \bibinfo{pages}{51--65}.
%Type = Article
\bibitem[{Bostan et~al.(2016)Bostan, Bousquet-M{\'e}lou, Kauers and Melczer}]{BBKM16}
\bibinfo{author}{Bostan, A.}, \bibinfo{author}{Bousquet-M{\'e}lou, M.}, \bibinfo{author}{Kauers, M.}, \bibinfo{author}{Melczer, S.}, \bibinfo{year}{2016}.
\newblock \bibinfo{title}{On 3-dimensional lattice walks confined to the positive octant}.
\newblock \bibinfo{journal}{Ann. Comb.} \bibinfo{volume}{20}, \bibinfo{pages}{661--704}.
\newblock \DOIprefix\doi{10.1007/s00026-016-0328-7}.
%Type = Article
\bibitem[{Bostan et~al.(2021)Bostan, Bousquet-M{\'e}lou and Melczer}]{BBMM21}
\bibinfo{author}{Bostan, A.}, \bibinfo{author}{Bousquet-M{\'e}lou, M.}, \bibinfo{author}{Melczer, S.}, \bibinfo{year}{2021}.
\newblock \bibinfo{title}{Counting walks with large steps in an orthant}.
\newblock \bibinfo{journal}{J. Eur. Math. Soc. (JEMS)} \bibinfo{volume}{23}, \bibinfo{pages}{2221--2297}.
\newblock \DOIprefix\doi{10.4171/JEMS/1053}.
%Type = Article
\bibitem[{Bostan et~al.(2017)Bostan, Chyzak, van Hoeij, Kauers and Pech}]{BCHKP}
\bibinfo{author}{Bostan, A.}, \bibinfo{author}{Chyzak, F.}, \bibinfo{author}{van Hoeij, M.}, \bibinfo{author}{Kauers, M.}, \bibinfo{author}{Pech, L.}, \bibinfo{year}{2017}.
\newblock \bibinfo{title}{Hypergeometric expressions for generating functions of walks with small steps in the quarter plane}.
\newblock \bibinfo{journal}{Eur. J. Comb.} \bibinfo{volume}{61}, \bibinfo{pages}{242--275}.
\newblock \DOIprefix\doi{doi:10.1016/j.ejc.2016.10.010}.
%Type = Article
\bibitem[{Bostan et~al.(2010)Bostan, Kauers and Van~Hoeij}]{BKvH10}
\bibinfo{author}{Bostan, A.}, \bibinfo{author}{Kauers, M.}, \bibinfo{author}{Van~Hoeij, M.}, \bibinfo{year}{2010}.
\newblock \bibinfo{title}{The complete generating function for {Gessel} walks is algebraic}.
\newblock \bibinfo{journal}{Proc. Am. Math. Soc.} \bibinfo{volume}{138}, \bibinfo{pages}{3063--3078}.
\newblock \DOIprefix\doi{10.1090/S0002-9939-2010-10398-2}.
%Type = Article
\bibitem[{Bostan et~al.(2014)Bostan, Raschel and Salvy}]{BRS14}
\bibinfo{author}{Bostan, A.}, \bibinfo{author}{Raschel, K.}, \bibinfo{author}{Salvy, B.}, \bibinfo{year}{2014}.
\newblock \bibinfo{title}{Non-{D}-finite excursions in the quarter plane}.
\newblock \bibinfo{journal}{J. Comb. Theory, Ser. A} \bibinfo{volume}{121}, \bibinfo{pages}{45--63}.
\newblock \DOIprefix\doi{10.1016/j.jcta.2013.09.005}.
%Type = Article
\bibitem[{Bousquet-M{\'e}lou(2005)}]{BM05}
\bibinfo{author}{Bousquet-M{\'e}lou, M.}, \bibinfo{year}{2005}.
\newblock \bibinfo{title}{Walks in the quarter plane: {Kreweras}' algebraic model}.
\newblock \bibinfo{journal}{Ann. Appl. Probab.} \bibinfo{volume}{15}, \bibinfo{pages}{1451--1491}.
\newblock \DOIprefix\doi{10.1214/105051605000000052}.
%Type = Article
\bibitem[{Bousquet-M{\'e}lou(2016a)}]{BM16_Gessel}
\bibinfo{author}{Bousquet-M{\'e}lou, M.}, \bibinfo{year}{2016}a.
\newblock \bibinfo{title}{An elementary solution of {Gessel}'s walks in the quadrant}.
\newblock \bibinfo{journal}{Adv. Math.} \bibinfo{volume}{303}, \bibinfo{pages}{1171--1189}.
\newblock \DOIprefix\doi{10.1016/j.aim.2016.08.038}.
%Type = Article
\bibitem[{Bousquet-M{\'e}lou(2016b)}]{BM16}
\bibinfo{author}{Bousquet-M{\'e}lou, M.}, \bibinfo{year}{2016}b.
\newblock \bibinfo{title}{Square lattice walks avoiding a quadrant}.
\newblock \bibinfo{journal}{J. Comb. Theory, Ser. A} \bibinfo{volume}{144}, \bibinfo{pages}{37--79}.
\newblock \DOIprefix\doi{10.1016/j.jcta.2016.06.010}.
%Type = Article
\bibitem[{Bousquet-M{\'e}lou(2023)}]{BM21}
\bibinfo{author}{Bousquet-M{\'e}lou, M.}, \bibinfo{year}{2023}.
\newblock \bibinfo{title}{Enumeration of three-quadrant walks via invariants: some diagonally symmetric models}.
\newblock \bibinfo{journal}{Can. J. Math.} \bibinfo{volume}{75}, \bibinfo{pages}{1566--1632}.
\newblock \DOIprefix\doi{10.4153/S0008414X22000487}.
%Type = Incollection
\bibitem[{Bousquet-M{\'e}lou and Mishna(2010)}]{MBM}
\bibinfo{author}{Bousquet-M{\'e}lou, M.}, \bibinfo{author}{Mishna, M.}, \bibinfo{year}{2010}.
\newblock \bibinfo{title}{Walks with small steps in the quarter plane}, in: \bibinfo{booktitle}{Algorithmic probability and combinatorics. Papers from the AMS special sessions, Chicago, IL, USA, October 5--6, 2007 and Vancouver, BC, Canada, October 4--5, 2008}. \bibinfo{publisher}{Providence, RI: American Mathematical Society (AMS)}, pp. \bibinfo{pages}{1--39}.
%Type = Article
\bibitem[{Buchacher et~al.(2021)Buchacher, Hofmanninger and Kauers}]{BHK21}
\bibinfo{author}{Buchacher, M.}, \bibinfo{author}{Hofmanninger, S.}, \bibinfo{author}{Kauers, M.}, \bibinfo{year}{2021}.
\newblock \bibinfo{title}{Walks with small steps in the 4d-orthant}.
\newblock \bibinfo{journal}{Ann. Comb.} \bibinfo{volume}{25}, \bibinfo{pages}{153--166}.
\newblock \DOIprefix\doi{10.1007/s00026-020-00520-5}.
%Type = Misc
\bibitem[{Chapon(2020)}]{Chapon_SW}
\bibinfo{author}{Chapon, F.}, \bibinfo{year}{2020}.
\newblock \bibinfo{howpublished}{Private Correspondence}.
%Type = Incollection
\bibitem[{Chapon et~al.(2020)Chapon, Fusy and Raschel}]{Poly}
\bibinfo{author}{Chapon, F.}, \bibinfo{author}{Fusy, {\'E}.}, \bibinfo{author}{Raschel, K.}, \bibinfo{year}{2020}.
\newblock \bibinfo{title}{Polyharmonic functions and random processes in cones}, in: \bibinfo{booktitle}{31st international conference on probabilistic, combinatorial and asymptotic methods for the analysis of algorithms, AofA 2020, Klagenfurt, Austria (virtual conference), June 15--19, 2020}. \bibinfo{publisher}{Wadern: Schloss Dagstuhl -- Leibniz Zentrum f{\"u}r Informatik}, p.~\bibinfo{pages}{19}.
\newblock \DOIprefix\doi{10.4230/LIPIcs.AofA.2020.9}. \bibinfo{note}{id/No 9}.
%Type = Article
\bibitem[{Connoll and Grupen(1993)}]{CG93}
\bibinfo{author}{Connoll, C.I.}, \bibinfo{author}{Grupen, R.A.}, \bibinfo{year}{1993}.
\newblock \bibinfo{title}{The applications of harmonic functions to robotics}.
\newblock \bibinfo{journal}{Journal of Field Robotics} \bibinfo{volume}{10}, \bibinfo{pages}{931--946}.
\newblock \DOIprefix\doi{https://doi.org/10.1002/rob.4620100704}.
%Type = Article
\bibitem[{Courtiel et~al.(2017)Courtiel, Melczer, Mishna and Raschel}]{CMMR17}
\bibinfo{author}{Courtiel, J.}, \bibinfo{author}{Melczer, S.}, \bibinfo{author}{Mishna, M.}, \bibinfo{author}{Raschel, K.}, \bibinfo{year}{2017}.
\newblock \bibinfo{title}{Weighted lattice walks and universality classes}.
\newblock \bibinfo{journal}{J. Comb. Theory, Ser. A} \bibinfo{volume}{152}, \bibinfo{pages}{255--302}.
\newblock \DOIprefix\doi{10.1016/j.jcta.2017.06.008}.
%Type = Misc
\bibitem[{Denisov et~al.(2024)Denisov, Tarasov and Wachtel}]{DTW23}
\bibinfo{author}{Denisov, D.}, \bibinfo{author}{Tarasov, A.}, \bibinfo{author}{Wachtel, V.}, \bibinfo{year}{2024}.
\newblock \bibinfo{title}{Expansions for random walks conditioned to stay positive}.
\newblock \href{http://arxiv.org/abs/2401.09929}{{\tt arXiv:2401.09929}}.
%Type = Article
\bibitem[{Denisov and Wachtel(2015)}]{Denisov}
\bibinfo{author}{Denisov, D.}, \bibinfo{author}{Wachtel, V.}, \bibinfo{year}{2015}.
\newblock \bibinfo{title}{Random walks in cones}.
\newblock \bibinfo{journal}{Ann. Probab.} \bibinfo{volume}{43}, \bibinfo{pages}{992--1044}.
\newblock \DOIprefix\doi{10.1214/13-AOP867}.
%Type = Article
\bibitem[{Dreyfus et~al.(2020)Dreyfus, Hardouin, Roques and Singer}]{SingerGenus0}
\bibinfo{author}{Dreyfus, T.}, \bibinfo{author}{Hardouin, C.}, \bibinfo{author}{Roques, J.}, \bibinfo{author}{Singer, M.F.}, \bibinfo{year}{2020}.
\newblock \bibinfo{title}{Walks in the quarter plane: genus zero case}.
\newblock \bibinfo{journal}{J. Comb. Theory, Ser. A} \bibinfo{volume}{174}, \bibinfo{pages}{24}.
\newblock \DOIprefix\doi{10.1016/j.jcta.2020.105251}. \bibinfo{note}{id/No 105251}.
%Type = Misc
\bibitem[{Elvey-Price et~al.(2023)Elvey-Price, Nessmann and Raschel}]{KilianAndrewPrivate}
\bibinfo{author}{Elvey-Price, A.}, \bibinfo{author}{Nessmann, A.}, \bibinfo{author}{Raschel, K.}, \bibinfo{year}{2023}.
\newblock \bibinfo{title}{Logarithmic terms in discrete heat kernel expansions in the quadrant}.
\newblock \href{http://arxiv.org/abs/2309.15209}{{\tt arXiv:2309.15209}}.
%Type = Article
\bibitem[{Fayolle and Raschel(2011)}]{FR11}
\bibinfo{author}{Fayolle, G.}, \bibinfo{author}{Raschel, K.}, \bibinfo{year}{2011}.
\newblock \bibinfo{title}{Random walks in the quarter-plane with zero drift: an explicit criterion for the finiteness of the associated group}.
\newblock \bibinfo{journal}{Markov Process. Relat. Fields} \bibinfo{volume}{17}, \bibinfo{pages}{619--636}.
%Type = Incollection
\bibitem[{Fayolle and Raschel(2012)}]{FR12}
\bibinfo{author}{Fayolle, G.}, \bibinfo{author}{Raschel, K.}, \bibinfo{year}{2012}.
\newblock \bibinfo{title}{Some exact asymptotics in the counting of walks in the quarter plane}, in: \bibinfo{booktitle}{Proceeding of the 23rd international meeting on probabilistic, combinatorial, and asymptotic methods in the analysis of algorithms (AofA'12), Montreal, Canada, June 18--22, 2012}. \bibinfo{publisher}{Nancy: The Association. Discrete Mathematics \& Theoretical Computer Science (DMTCS)}, pp. \bibinfo{pages}{109--124}.
%Type = Article
\bibitem[{Fayolle and Raschel(2015a)}]{FR17}
\bibinfo{author}{Fayolle, G.}, \bibinfo{author}{Raschel, K.}, \bibinfo{year}{2015}a.
\newblock \bibinfo{title}{About a possible analytic approach for walks in the quarter plane with arbitrary big jumps}.
\newblock \bibinfo{journal}{C. R., Math., Acad. Sci. Paris} \bibinfo{volume}{353}, \bibinfo{pages}{89--94}.
\newblock \DOIprefix\doi{10.1016/j.crma.2014.11.015}.
%Type = Article
\bibitem[{Fayolle and Raschel(2015b)}]{FR15}
\bibinfo{author}{Fayolle, G.}, \bibinfo{author}{Raschel, K.}, \bibinfo{year}{2015}b.
\newblock \bibinfo{title}{About a possible analytic approach for walks in the quarter plane with arbitrary big jumps}.
\newblock \bibinfo{journal}{C. R., Math., Acad. Sci. Paris} \bibinfo{volume}{353}, \bibinfo{pages}{89--94}.
\newblock \URLprefix \url{hal.inria.fr/hal-01021327/file/1-s2.0-S1631073X14003082-main.pdf}, \DOIprefix\doi{10.1016/j.crma.2014.11.015}.
%Type = Article
\bibitem[{Feierl(2014)}]{F14}
\bibinfo{author}{Feierl, T.}, \bibinfo{year}{2014}.
\newblock \bibinfo{title}{Asymptotics for the number of walks in a {Weyl} chamber of type {{\(B\)}}}.
\newblock \bibinfo{journal}{Random Struct. Algorithms} \bibinfo{volume}{45}, \bibinfo{pages}{261--305}.
\newblock \DOIprefix\doi{10.1002/rsa.20467}.
%Type = Book
\bibitem[{Flajolet and Sedgewick(2009)}]{Flajolet}
\bibinfo{author}{Flajolet, P.}, \bibinfo{author}{Sedgewick, R.}, \bibinfo{year}{2009}.
\newblock \bibinfo{title}{Analytic Combinatorics}.
\newblock \bibinfo{publisher}{Cambridge University Press}.
%Type = Article
\bibitem[{Flatto and Hahn(1984)}]{FH84}
\bibinfo{author}{Flatto, L.}, \bibinfo{author}{Hahn, S.}, \bibinfo{year}{1984}.
\newblock \bibinfo{title}{Two parallel queues created by arrivals with two demands. {I}}.
\newblock \bibinfo{journal}{SIAM J. Appl. Math.} \bibinfo{volume}{44}, \bibinfo{pages}{1041--1053}.
\newblock \DOIprefix\doi{10.1137/0144074}.
%Type = Article
\bibitem[{Hardouin and Singer(2021)}]{Singer}
\bibinfo{author}{Hardouin, C.}, \bibinfo{author}{Singer, M.F.}, \bibinfo{year}{2021}.
\newblock \bibinfo{title}{On differentially algebraic generating series for walks in the quarter plane}.
\newblock \bibinfo{journal}{Sel. Math., New Ser.} \bibinfo{volume}{27}, \bibinfo{pages}{49}.
\newblock \DOIprefix\doi{10.1007/s00029-021-00703-9}. \bibinfo{note}{id/No 89}.
%Type = Article
\bibitem[{Hoang et~al.(2022)Hoang, Raschel and Tarrago}]{Hung}
\bibinfo{author}{Hoang, V.H.}, \bibinfo{author}{Raschel, K.}, \bibinfo{author}{Tarrago, P.}, \bibinfo{year}{2022}.
\newblock \bibinfo{title}{Constructing discrete harmonic functions in wedges}.
\newblock \bibinfo{journal}{Trans. Am. Math. Soc.} \bibinfo{volume}{375}, \bibinfo{pages}{4741--4782}.
\newblock \DOIprefix\doi{10.1090/tran/8615}.
%Type = Incollection
\bibitem[{Kauers and Yatchak(2015)}]{KY15}
\bibinfo{author}{Kauers, M.}, \bibinfo{author}{Yatchak, R.}, \bibinfo{year}{2015}.
\newblock \bibinfo{title}{Walks in the quarter plane with multiple steps}, in: \bibinfo{booktitle}{Proceedings of the 27th international conference on formal power series and algebraic combinatorics, FPSAC 2015, Daejeon, South Korea, July 6--10, 2015}. \bibinfo{publisher}{Nancy: The Association. Discrete Mathematics \& Theoretical Computer Science (DMTCS)}, pp. \bibinfo{pages}{25--36}.
%Type = Article
\bibitem[{Kurkova and Raschel(2012)}]{KR12}
\bibinfo{author}{Kurkova, I.}, \bibinfo{author}{Raschel, K.}, \bibinfo{year}{2012}.
\newblock \bibinfo{title}{On the functions counting walks with small steps in the quarter plane}.
\newblock \bibinfo{journal}{Publ. Math., Inst. Hautes {\'E}tud. Sci.} \bibinfo{volume}{116}, \bibinfo{pages}{69--114}.
\newblock \DOIprefix\doi{10.1007/s10240-012-0045-7}.
%Type = Book
\bibitem[{Lurie and Vasiliev(1995)}]{Elasticity}
\bibinfo{author}{Lurie, S.A.}, \bibinfo{author}{Vasiliev, V.V.}, \bibinfo{year}{1995}.
\newblock \bibinfo{title}{The Biharmonic Problem in the Theory of Elasticity}.
\newblock \bibinfo{publisher}{Taylor \& Francis Ltd}.
%Type = Article
\bibitem[{Maschke(1898)}]{Mas98}
\bibinfo{author}{Maschke, H.}, \bibinfo{year}{1898}.
\newblock \bibinfo{title}{Ueber den arithmetischen charakter der coefficienten der substitutionen endlicher linearer substitutionsgruppen.}
\newblock \bibinfo{journal}{Math. Ann.} \bibinfo{volume}{50}, \bibinfo{pages}{492--498}.
\newblock \DOIprefix\doi{10.1007/BF01444297}.
%Type = Book
\bibitem[{Melczer(2021)}]{ACSV}
\bibinfo{author}{Melczer, S.}, \bibinfo{year}{2021}.
\newblock \bibinfo{title}{An Invitation to Analytic Combinatorics}.
\newblock \bibinfo{publisher}{Springer International Publishing AG}.
%Type = Article
\bibitem[{Melczer and Mishna(2016)}]{MM16}
\bibinfo{author}{Melczer, S.}, \bibinfo{author}{Mishna, M.}, \bibinfo{year}{2016}.
\newblock \bibinfo{title}{Asymptotic lattice path enumeration using diagonals}.
\newblock \bibinfo{journal}{Algorithmica} \bibinfo{volume}{75}, \bibinfo{pages}{782--811}.
\newblock \DOIprefix\doi{10.1007/s00453-015-0063-1}.
%Type = Article
\bibitem[{Melczer and Wilson(2019)}]{MW19}
\bibinfo{author}{Melczer, S.}, \bibinfo{author}{Wilson, M.C.}, \bibinfo{year}{2019}.
\newblock \bibinfo{title}{Higher dimensional lattice walks: connecting combinatorial and analytic behavior}.
\newblock \bibinfo{journal}{SIAM J. Discrete Math.} \bibinfo{volume}{33}, \bibinfo{pages}{2140--2174}.
\newblock \DOIprefix\doi{10.1137/18M1220856}.
%Type = Article
\bibitem[{Mishna and Rechnitzer(2009)}]{MR09}
\bibinfo{author}{Mishna, M.}, \bibinfo{author}{Rechnitzer, A.}, \bibinfo{year}{2009}.
\newblock \bibinfo{title}{Two non-holonomic lattice walks in the quarter plane}.
\newblock \bibinfo{journal}{Theor. Comput. Sci.} \bibinfo{volume}{410}, \bibinfo{pages}{3616--3630}.
\newblock \DOIprefix\doi{10.1016/j.tcs.2009.04.008}.
%Type = Misc
\bibitem[{Nessmann(2022)}]{Nes22}
\bibinfo{author}{Nessmann, A.}, \bibinfo{year}{2022}.
\newblock \bibinfo{title}{Polyharmonic functions in the quarter plane}.
\newblock \href{http://arxiv.org/abs/2212.07258}{{\tt arXiv:2212.07258}}.
%Type = Misc
\bibitem[{Nessmann(2023)}]{Nes23}
\bibinfo{author}{Nessmann, A.}, \bibinfo{year}{2023}.
\newblock \bibinfo{title}{Full asymptotic expansion for orbit-summable quadrant walks and discrete polyharmonic functions}.
\newblock \href{http://arxiv.org/abs/2307.11539}{{\tt arXiv:2307.11539}}.
%Type = Article
\bibitem[{Price(2020)}]{EP20}
\bibinfo{author}{Price, A.E.}, \bibinfo{year}{2020}.
\newblock \bibinfo{title}{Counting lattice walks by winding angle}.
\newblock \bibinfo{journal}{S{\'e}min. Lothar. Comb.} \bibinfo{volume}{84B}, \bibinfo{pages}{12}.
\newblock \URLprefix \url{www.mat.univie.ac.at/~slc/wpapers/FPSAC2020//43.html}. \bibinfo{note}{id/No 43}.
%Type = Article
\bibitem[{Raschel(2011)}]{Ras11}
\bibinfo{author}{Raschel, K.}, \bibinfo{year}{2011}.
\newblock \bibinfo{title}{Green functions for killed random walks in the {Weyl} chamber of {{\(\mathrm{Sp}(4)\)}}}.
\newblock \bibinfo{journal}{Ann. Inst. Henri Poincar{\'e}, Probab. Stat.} \bibinfo{volume}{47}, \bibinfo{pages}{1001--1019}.
\newblock \DOIprefix\doi{10.1214/10-AIHP405}.
%Type = Article
\bibitem[{Raschel(2012)}]{Kilian_JEMS}
\bibinfo{author}{Raschel, K.}, \bibinfo{year}{2012}.
\newblock \bibinfo{title}{Counting walks in a quadrant: a unified approach via boundary value problems}.
\newblock \bibinfo{journal}{J. Eur. Math. Soc. (JEMS)} \bibinfo{volume}{14}, \bibinfo{pages}{749--777}.
\newblock \DOIprefix\doi{10.4171/JEMS/317}.
%Type = Article
\bibitem[{Raschel(2014)}]{Conformal}
\bibinfo{author}{Raschel, K.}, \bibinfo{year}{2014}.
\newblock \bibinfo{title}{Random walks in the quarter plane, discrete harmonic functions and conformal mappings}.
\newblock \bibinfo{journal}{Stochastic Processes Appl.} \bibinfo{volume}{124}, \bibinfo{pages}{3147--3178}.
\newblock \DOIprefix\doi{10.1016/j.spa.2014.04.013}.
%Type = Article
\bibitem[{Raschel and Trotignon(2019)}]{RT19}
\bibinfo{author}{Raschel, K.}, \bibinfo{author}{Trotignon, A.}, \bibinfo{year}{2019}.
\newblock \bibinfo{title}{On walks avoiding a quadrant}.
\newblock \bibinfo{journal}{Electron. J. Comb.} \bibinfo{volume}{26}, \bibinfo{pages}{research paper p3.31, 34}.
\newblock \URLprefix \url{www.combinatorics.org/ojs/index.php/eljc/article/view/v26i3p31}.
%Type = Article
\bibitem[{Trotignon(2022)}]{Tro22}
\bibinfo{author}{Trotignon, A.}, \bibinfo{year}{2022}.
\newblock \bibinfo{title}{Discrete harmonic functions in the three-quarter plane}.
\newblock \bibinfo{journal}{Potential Anal.} \bibinfo{volume}{56}, \bibinfo{pages}{267--296}.
\newblock \DOIprefix\doi{10.1007/s11118-020-09884-y}.
%Type = Incollection
\bibitem[{Yatchak(2017)}]{Yat17}
\bibinfo{author}{Yatchak, R.}, \bibinfo{year}{2017}.
\newblock \bibinfo{title}{Automated positive part extraction for lattice path generating functions in the octant}, in: \bibinfo{booktitle}{Extended abstracts of the ninth European conference on combinatorics, graph theory and applications, EuroComb 2017, Vienna, Austria, August 28 -- September 1, 2017}. \bibinfo{publisher}{Amsterdam: Elsevier}, pp. \bibinfo{pages}{1061--1067}.
\newblock \DOIprefix\doi{10.1016/j.endm.2017.07.073}.

\end{thebibliography}

\appendix

\section{An example with large steps}\label{app:largesteps}
Consider the model with steps $(1,0),(-1,0),(0,-1),(-2,1)$. In~\cite[Prop.~16]{BBMM21} it is shown that \begin{align*}
    Q(x,y;t)=\left[x^{>}y^{>}\right]\frac{(x^2 + 1) (x + y) (y - x) (x^2y - 2 x - y) (x^3 - x - 2 y)}{x^7y^3(1-tS(x,y))},
\end{align*}
with $S(x,y)$ the step counting polynomial \begin{align*}
    S(x,y)=x+x^{-1}+y^{-1}+x^{-2}y.
\end{align*}
We find the dominant saddle point to be at $s_0:=\left(\sqrt{3},\sqrt{3}\right)$, and one other saddle point at $s_1=\left(-\sqrt{3},-\sqrt{3}\right)$ associated to it. We can check that we have\footnote{This is in fact a consequence of $K(x,y)$ still being quadratic in $y$, which allows one to argue in a fashion similar as in the proof of Thm.~\ref{thm:mainthm}.}, using the notation as in the proof of Thm.~\ref{thm:mainthm}, $N(s_0)=0$ (in particular it is not infinite), thus we can proceed in the same manner as in aforementioned theorem. We obtain \begin{align*}
    \gamma=2\sqrt{3},\quad c=3,
\end{align*}
and have \begin{small}\begin{align*}
    v_1(k,l)=&\frac{16}{\pi}\cdot \sqrt{3}^{-1 - k - l} (1 + k) (1 + l) (3 + k + 2 l),\\
    v_2(k,l)=&-\frac{2}{\pi} \sqrt{3}^{-1 - k - l}(1 + k) (1 + l) (3 + k + 2 l) (107 + 4 k^2 + 32 l + 16 l^2 + 
   8 k (1 + l)) ,\\
   v_3(k,l)=&\frac{1}{8\pi} \sqrt{3}^{-3 - k - l} (1 + k) (1 + l) (3 + k + 2 l) (15205 + 16 k^4 + 8672 l + 4976 l^2 + 
   832 l^3 \\&+ 256 l^4 + 64 k^3 (1 + l) + 8 k^2 (157 + 48 l + 24 l^2) + 
   16 k (149 + 157 l + 36 l^2 + 16 l^3)) .
\end{align*}\end{small}
\section{A three-dimensional example}\label{app:3dim}
Consider the model with steps $(-1,-1,-1), (-1,-1,1), (-1,1,0),(1,0,0)$. In~\cite[4.3]{BBKM16} it is shown that we have \begin{footnotesize}\begin{equation*}
Q(x,y,z;t)=\left[x^{>0}\right]\left[y^{>0}\right]\left[z^{>0}\right]\frac{\left(x-x^{-1}y-x^{-1}y^{-1}z-x^{-1}y^{-1}z^{-1}\right)\left(y-y^{-1}z-y^{-1}z^{-1}\right)\left(z-z^{-1}\right)}{xyz\left(1-tS(x,y,z)\right)},
\end{equation*}
\end{footnotesize}
with the step counting polynomial \begin{align*}
S(x,y,z)=x^{-1}y^{-1}z^{-1}+x^{-1}y^{-1}z+x^{-1}y+x.
\end{align*}
We can find the dominant saddle point of $S(x,y,z)$ to be at $s_0=\left(2^{3/4},2^{1/2},1\right)$, with $7$ others associated to it. Using the notation as in (\ref{eq:mainthm:asymptotics}), we can check that \begin{align*}
\gamma=2\cdot 2^{3/4},\quad c =\frac{7}{2},
\end{align*}
and that we have \begin{small} \begin{align*}
v_1(k,l,m)=&\frac{2^6}{\pi^{3/2}}2^{ - 3k/4 - l/2} (1 + k) (1 + l) (1 + m),\\
v_2(k,l,m)=&-\frac{2^4}{\pi^{3/2}}2^{- 3k/4 - l/2} (1 + k) (1 + l) (1 + m) (63 - 8 k + 2 k^2 - 
   4 l + 4 l^2 + 16 m + 8 m^2),\\
   v_3(k,l,m)=&\frac{2}{\pi^{3/2}}2^{ - 3k/4 - l/
  2} (1 + k) (1 + l) (1 + m) (5313 - 32 k^3 + 4 k^4 - 32 l^3 + 
   16 l^4 + 3040 m \\&+ 1776 m^2 + 256 m^3 + 64 m^4 - 
   32 k (43 - 3 l + 3 l^2 + 12 m + 6 m^2) \\&+ 
   8 l^2 (93 + 16 m + 8 m^2) + 
   4 k^2 (99 - 4 l + 4 l^2 + 16 m + 8 m^2) - 
   8 l (103 + 48 m + 24 m^2)).
\end{align*}\end{small}

\end{document}